\documentclass[11pt]{article}

\usepackage[margin=1.25in]{geometry}
\usepackage{amsmath,amssymb,amsthm,mathtools,bm}
\usepackage{algorithm}
\usepackage{algorithmic}
\usepackage{microtype}
\usepackage{graphicx}
\usepackage{xcolor}
\usepackage{caption}
\usepackage[numbers,square,sort&compress]{natbib}
\usepackage[hidelinks]{hyperref}
\usepackage{xurl}
\usepackage{subcaption}
\usepackage{comment}
\usepackage{placeins}
\usepackage{booktabs}

\newtheorem{assumption}{Assumption}
\newtheorem{lemma}{Lemma}
\newtheorem{theorem}{Theorem}
\newtheorem{corollary}{Corollary}

\newtheorem{remark}{Remark}

\newcommand{\R}{\mathbb R}
\newcommand{\E}{\mathbb E}
\newcommand{\F}{\mathcal F}
\newcommand{\cE}{\mathcal E}
\newcommand{\diag}{\operatorname{diag}}
\newcommand{\Tr}{\operatorname{Tr}}
\newcommand{\email}[1]{\texttt{#1}}

\title{Variance-Reduced Stochastic Hessian-Driven Nesterov Acceleration with Adam-Style Preconditioning}
\author{Yaxin Yu\thanks{School of Mathematics, Sichuan University, Chengdu, China (\email{yaxin\_yu@stu.scu.edu.cn}, \email{fmf@scu.edu.cn}).}
\and Long Chen\thanks{Corresponding author. Department of Mathematics, University of California, Irvine, CA, USA (\email{lchen7@uci.edu}).}
\and Minfu Feng\footnotemark[1]}
\date{}

\begin{document}
	\maketitle

\begin{abstract}
	Stochastic first-order methods use sampled gradients to reduce
	per-iteration costs, but estimation errors can slow convergence
	and, in some settings, cause divergence. To improve gradient accuracy
	while retaining acceleration, we develop scalar and diagonal versions
	of {\it VR-Adam-SHANG} for smooth stochastic convex optimization.
	The methods combine STORM-type recursive estimation, Adam-style
	preconditioning, and SHANG's accelerated interpolation and correction,
	using at most two stochastic gradient evaluations per iteration
	without full-gradient refreshes. Under unbiased sampling, bounded
	variance, and mean-square smoothness of gradient noise, polynomial
	steps with prescribed initialization yield $O(1/K)$ mean-square
	gradient-estimation error and expected objective-gap bounds combining
	$O(1/K^2)$ initial-energy decay with an $O(1/\sqrt K)$ stochastic
	term after $K$ iterations. These guarantees require neither monotonicity of the preconditioner
	nor a uniform positive lower bound on its eigenvalues.
	The analysis also establishes polynomial-step convergence of
	projected nonrecursive {\it Adam-SHANG-s}, filling a gap in its
	original theory. Experiments compare the methods at matched sample-gradient budgets.
	In every tested setting, at least one of the two versions achieves
	the lowest mean terminal objective gap among all compared methods.
\end{abstract}

\section{Introduction}\label{sec:introduction}
Many optimization problems involve minimizing an expected objective
whose gradient is costly to evaluate exactly. They can often be formulated as
\begin{equation}\label{eq:problem}
	\min_{x\in\R^d} f(x), \qquad f(x) :=\E_{\xi\sim\mathcal D}[f(x,\xi)].
\end{equation}
Here the objective $f$ is convex and smooth, and we assume
that a minimizer $x^\star$ exists. Examples include regularized
logistic regression and convex learning from streaming data.
Stochastic gradient descent (SGD) replaces the exact gradient by a
sampled gradient \citep{robbins1951stochastic}:
\[
x_{k+1}=x_k-\eta_k\nabla f(x_k,\xi_k).
\]
When the gradient-noise variance does not vanish near a minimizer,
constant steps can leave a nonzero objective gap. Decreasing steps
reduces the effect of noise but slows progress, while larger batches
increase the cost per iteration. Uneven curvature and gradient scales
pose further difficulties. Acceleration, adaptive preconditioning,
and variance reduction address these difficulties through momentum,
gradient-dependent scaling, and more accurate gradient estimates,
respectively. This paper combines variance reduction with acceleration
and adaptive preconditioning, using recursive gradient estimates to
drive both the accelerated updates and the preconditioner.

\paragraph{Acceleration and noise stability.}
For smooth convex objectives with exact gradients, Nesterov acceleration \citep{Nesterov1983} improves the worst-case objective-gap rate of gradient descent from $O(1/K)$ to $O(1/K^2)$ after $K$ iterations. Simply replacing exact gradients by stochastic estimates in NAG need
not preserve acceleration \citep{kidambi2018insufficiency,LiuBelkin2020}; standard Nesterov SGD can even diverge at step sizes for which SGD converges \citep{LiuBelkin2020}. Noise-aware accelerated methods achieve an $O(1/K^2)$ deterministic term and an $O(1/\sqrt K)$ stochastic term under bounded variance \citep{hu2009accelerated,lan2012acsa}. Further acceleration results exploit least-squares structure \citep{jain2018} or multiplicative noise \citep{VaswaniBachSchmidt2019,GuptaSiegelWojtowytsch2024, hermant2025gradientcorrelation}.

The Hessian-driven Nesterov accelerated gradient (HNAG) flow \citep{ChenLuo2019optimization} provides a Lyapunov framework for acceleration. SHANG and SHANG++ \citep{yu2026shangrobuststochasticacceleration} use this structure to improve stability under multiplicative noise. SHANG uses a semi-implicit discretization; SHANG++ adds damping. Their accelerated guarantees use a multiplicative-noise scaling assumption, under which the noise variance is bounded by a multiple of $\|\nabla f(x_k)\|^2$. We retain the accelerated update structure while allowing bounded-variance noise that need not vanish at a minimizer.

\paragraph{Adaptive scaling.}
Adaptive methods use gradient history to accommodate heterogeneous
gradient scales and feature frequencies. AdaGrad \citep{duchi2011adaptive} accumulates squared gradients, whereas Adam \citep{kingma2015adam} and AMSGrad \citep{reddi2018convergence} combine momentum with adaptive scaling based on exponentially averaged squared gradients. Shampoo \citep{gupta2018shampoo} goes beyond diagonal scaling
with structured matrix preconditioners. AcceleGrad \citep{levy2018online} accelerates with adaptive scalar
steps without knowing the smoothness constant, while diagonal
adaptation is combined with acceleration for stochastic convex
optimization in~\citep{deng2018optimaladaptive,ene2021adaptive}. 
Adam is a central motivation for this work: it combines momentum
with coordinatewise scaling, allowing updates to respond to different
gradient magnitudes across coordinates. Understanding how to retain
this adaptation while proving convergence motivates the
Adam-type methods below.

{\it Adam-HNAG} \citep{yu2026adamhnag} couples the HNAG flow with
feedback based on squared gradients: gradient magnitudes drive the
preconditioner update, and the inverse preconditioner scales the
gradient response. This feedback balances changes in the metric
against energy descent in the Lyapunov analysis.
{\it Adam-SHANG} \citep{yu2026adamshangconvergentadamtypemethod} extends this approach to stochastic gradients. Its synchronous variant, {\it Adam-SHANG-s}, uses the current gradient in both the direction and the preconditioner and thus is closer to Adam. The original analysis in  \citep{yu2026adamshangconvergentadamtypemethod}, however, treats a lagged preconditioner and leaves convergence of the synchronous variant open.

\paragraph{Variance reduction.}
Adaptive preconditioning changes how gradient estimates enter the
updates, while variance reduction improves the estimates themselves
by exploiting past information and correlated samples.
For finite sums, SVRG \citep{johnson2013accelerating} uses a reference point and its full gradient, while SAGA \citep{defazio2014saga} stores component gradients. SARAH \citep{nguyen2017sarah} and SPIDER \citep{fang2018spider} update estimates through gradient differences, with periodic full-gradient or large-batch refreshes.
STORM \citep{cutkosky2019momentum} replaces these refreshes
by a damped recursion with $0\le\rho_k<1$:
$$
g_{k+1}=\nabla f(x_{k+1},\xi_{k+1})
+\rho_k\bigl(g_k-\nabla f(x_k,\xi_{k+1})\bigr).
$$
Here $g_k$ estimates $\nabla f(x_k)$, and both gradients use the
\emph{same} fresh sample or minibatch. The recursion reflects
the bias--variance tradeoff in statistical learning
\citep[Section~7.3]{hastie2009elements}: a biased estimate can have
smaller mean-square error if the variance reduction outweighs the
squared bias. Here the decomposition is conditional on the past.
Let $U_k=g_k-\nabla f(x_k)$ denote the gradient-estimation error.
Unbiased sampling gives $\E_kU_{k+1}=\rho_kU_k$, where $\E_k$
conditions on the history before drawing $\xi_{k+1}$.
The conditional mean-square error therefore decomposes as
\begin{equation}\label{intro:tradeoff}
\begin{aligned}
\E_k\|U_{k+1}\|^2
&=\underbrace{\rho_k^2\|U_k\|^2}_{\text{squared conditional bias}}
+\underbrace{\E_k\|U_{k+1}-\rho_kU_k\|^2}_{\operatorname{Var}_k(U_{k+1})}.
\end{aligned}
\end{equation}
Increasing $\rho_k$ retains more past error but damps fresh noise
through the factor $1-\rho_k$. Under mean-square smoothness of
the gradient noise, the conditional second moment of the same-sample
noise difference is at most $O(\|x_{k+1}-x_k\|^2)$.
Thus small query motion---the distance between successive points
where gradients are evaluated---and damped fresh noise together
control the variance contribution.
Neither gradient tables nor periodic refreshes are required.

Variance reduction has already been combined with acceleration in Katyusha \citep{allenzhu2018katyusha}, Varag \citep{lan2019varag}, and related methods \citep{driggs2022accelerating,liu2022adaptivevr}. SUPER-ADAM \citep{huang2021superadam}, VRAdam \citep{li2023adam}, and MARS \citep{yuan2025mars} combine recursive estimation with adaptive preconditioning and prove nonconvex stationarity guarantees.
The closest convex comparison is $\mu^2$-ExtraSGD
\citep{dahan2025mu2sgd}, which combines recursive estimation,
weighted averaging, and Euclidean extragradient steps to obtain
the same objective-bound orders as in \eqref{eq:order}. It requires three stochastic gradient
evaluations per iteration and does not incorporate adaptive scaling.
Both versions of our method achieve these orders with at most
two evaluations per iteration. At matched sample-gradient budgets,
both attain lower mean terminal objective gaps than $\mu^2$-ExtraSGD
in our experiments.

\paragraph{Contributions and scope.}
We develop scalar and diagonal versions of {\it VR-Adam-SHANG},
combining STORM-type recursive estimation, Adam-style preconditioning,
and SHANG's accelerated interpolation and correction.
The main analytical challenge is a coupled recursion: the conditionally
biased gradient estimate determines its own preconditioner and drives
the query motion, which in turn governs subsequent estimation errors. 
\begin{enumerate}
\item We give a Lyapunov analysis under unbiased sampling, bounded
variance, and mean-square smoothness of gradient noise. 
Compared with SHANG, it controls the dependence between the estimate
and its preconditioner; compared with nonrecursive {\it Adam-SHANG-s},
it additionally handles the coupling between query motion and
recursive estimation error. Gradient feedback in the preconditioner
update bounds correction lengths, while projection bounds the
auxiliary variable, jointly controlling query motion and the gradient-estimation error. Closing these coupled bounds is the main technical contribution.

\item With the prescribed initialization and polynomial steps, both versions
achieve mean-square gradient-estimation error
$\E\| g_K - \nabla f(x_K)\|^2=O(1/K)$.
For the corrected output of the scalar version, we obtain
\begin{equation}\label{eq:order}
\E[f(x_K^+)-f(x^\star)]
=O\!\left(
\frac{\E[\cE(z_0^+,P_0)]}{K^2}
+\frac{R\sigma+L_\delta R^2}{\sqrt K}
\right).
\end{equation}
Here $\cE(z_0^+,P_0)$ is the initial energy, while $R$, $\sigma$,
and $L_\delta$ quantify the projection radius, noise magnitude,
and query dependence, respectively.
The diagonal version attains the same orders, with stochastic coefficient
$\sqrt d\,R\sigma+dL_\delta R^2$.
Query-independent additive noise gives $L_\delta=0$.

Under the same projection and polynomial steps, we also establish convergence of nonrecursive {\it Adam-SHANG-s} with the same stochastic order, addressing a gap in its convergence
theory in \cite{yu2026adamshangconvergentadamtypemethod}. 

\item Experiments on graph-signal denoising and logistic regression with synthetic and real data show reduced mean-square gradient-estimation errors consistent with the theory. In every tested setting, either the scalar or the diagonal version achieves the lowest mean terminal objective gap at matched sample-gradient budgets, outperforming SHANG, {\it Adam-SHANG-s}, STORM, $\mu^2$-ExtraSGD, and Adam. 
\end{enumerate}

\paragraph{Notation and assumptions.}
Unsubscripted norms are Euclidean. A symmetric matrix $P\succ0$ means $P$ is symmetric and positive definite. Write
$\langle x,y\rangle_P:=\langle Px,y\rangle$ and
$\|x\|_P^2:=\langle Px,x\rangle$.
We refer to $P_k$ as the preconditioner and apply its inverse
$P_k^{-1}$ to the gradient estimate.
We write $\diag(u.^2):=\diag(u_1^2,\ldots,u_d^2)$.
The Bregman divergence measures the error of linearization:
\begin{equation}\label{bregman}
	D_f(y,x):=f(y)-f(x)-\langle\nabla f(x),y-x\rangle.
\end{equation}
We impose convexity and an $L$-Lipschitz gradient on the mean objective:
\begin{equation}\label{ass:convexity}
	D_f(x,y)\ge0,\qquad
	\|\nabla f(x)-\nabla f(y)\|\le L\|x-y\|,
	\qquad x,y\in\R^d.
\end{equation}

Let $\F_k$ be the sigma-field generated by initialization and the
history through step $k$. Draw $\xi_{k+1}\sim\mathcal D$
independently of $\F_k$ and write
$\E_k[\cdot]:=\E[\cdot\mid\F_k]$.
For any square-integrable random vector $X$, define its conditional
total variance by
$\operatorname{Var}_k(X):=\E_k\|X-\E_kX\|^2$.
The first two assumptions hold for every $k\ge0$ and
$\F_k$-measurable $x\in\R^d$.

\begin{assumption}[Unbiasedness]\label{ass:unbiased-oracle}
$\E_k[\nabla f(x,\xi_{k+1})]=\nabla f(x)$.
\end{assumption}

\begin{assumption}[Coordinatewise bounded variance]
	\label{ass:bounded-variance}
	There are $\sigma_i\ge0$, with
	\mbox{$\sigma^2:=\sum_{i=1}^d\sigma_i^2$}, such that
	\begin{equation}\label{eq:metric-variance}
		\E_k[(\nabla f(x,\xi_{k+1})-\nabla f(x))_i^2]
		\le\sigma_i^2,\qquad i=1,\ldots,d.
	\end{equation}
\end{assumption}

\begin{assumption}[Mean-square smoothness of gradient noise]
	\label{ass:smoothness}
	There is $L_\delta\ge0$ such that, for all $x,y\in\R^d$,
	\begin{equation}\label{eq:centered-smoothness-constant}
		\E_\xi\bigl\|
		\nabla f(x,\xi)-\nabla f(y,\xi)
		-(\nabla f(x)-\nabla f(y))
		\bigr\|^2
		\le L_\delta^2\|x-y\|^2.
	\end{equation}
\end{assumption}
Query-independent additive noise and exact gradients give $L_\delta=0$.

\section{Scalar VR-Adam-SHANG}
We first present the scalar version of {\it VR-Adam-SHANG},
combining variance-reduced Hessian-driven Nesterov acceleration
with scalar Adam-style preconditioning. This setting clarifies
the main argument; Section~\ref{sec:vr-adam-shang} extends the
method to coordinatewise scaling.

\subsection{Algorithm}
Choose $R>0$ such that $\|x^\star\|\le R/2$, and let $\mathcal B_R=\{y\in\R^d:\|y\|\le R/2\}$. Initialize $x_0,x_0^+,y_0\in\mathcal B_R$, $p_0>0$ and $\rho_0=0$. For $k \ge 0$, the updates are
\begin{subequations}\label{eq:algorithmvrshang}
\allowdisplaybreaks[4]
\begin{align}
\frac{x_{k+1}-x_k^+}{\alpha_k}
&=y_k-x_{k+1},\label{eq:algorithmvrshang-x}\\
g_{k+1}
&=\nabla f(x_{k+1},\xi_{k+1})
+\rho_k\bigl(g_k-\nabla f(x_k,\xi_{k+1})\bigr),
\qquad k\ge1,\label{eq:algorithmvrshang-g}\\
\frac{p_{k+1}-p_k}{\tilde\alpha_k}
&=-p_k+\gamma_kp_{k+1}^{-1}\|g_{k+1}\|^2,
\label{eq:algorithmvrshang-P}\\
\frac{y_{k+\frac12}-y_k}{\tilde\alpha_k}
&=-p_{k+1}^{-1}g_{k+1},
\label{eq:algorithmvrshang-y-half}\\
y_{k+1}
&=\operatorname*{argmin}_{y\in\mathcal B_R}
\|y-y_{k+\frac12}\|^2,\label{eq:algorithmvrshang-y}\\
x_{k+1}^+
&=x_{k+1}-\eta_{k+1}p_{k+1}^{-1}g_{k+1}.
\label{eq:algorithmvrshang-xplus}
\end{align}
\end{subequations}
Here $\alpha_k,\eta_{k+1}>0$ are the discretization and correction
steps, respectively; $\tilde\alpha_k=\alpha_k/(1+\alpha_k)\in(0,1)$, $\gamma_k\ge0$ controls gradient feedback,
and $\rho_k\in[0,1)$ is the memory weight. At $k=0$, set $g_1=\nabla f(x_1,\xi_1)$; no $g_0$ is needed. The preconditioner update \eqref{eq:algorithmvrshang-P}
preserves positivity and admits the explicit solution given in
Algorithm~\ref{alg:vr-shang}.

\begin{algorithm}[H]
\caption{{\it VR-Adam-SHANG}: scalar version}
\label{alg:vr-shang}
\renewcommand{\baselinestretch}{1.12}
\small
\begin{algorithmic}[1]
\REQUIRE $L>0$, $R>0$, $x_0,x_0^+,y_0\in\mathcal B_R$,
and iteration count $K\ge1$.
\par
\STATE Initialize $p_0\gets2L$ and $\rho_0\gets0$.
\FOR{$k=0,\ldots,K-1$}
\STATE $\tilde\alpha_k\gets2/(k+4)$.
\STATE $\alpha_k\gets\tilde\alpha_k/(1-\tilde\alpha_k)$,
$\gamma_k\gets\tilde\alpha_k/R^2$, $\eta_{k+1}\gets2\tilde\alpha_k^2$. \textbf{If} $k\ge1$, set $\rho_k\gets1-\tilde\alpha_k$.
\STATE $x_{k+1}\gets(x_k^++\alpha_k y_k)/(1+\alpha_k)$.
\STATE Draw a fresh sample $\xi_{k+1}\sim\mathcal D$, independently of $\F_k$,
and evaluate $\nabla f(x_{k+1},\xi_{k+1})$.
\STATE \textbf{If} $\rho_k=0$, set $g_{k+1}\gets\nabla f(x_{k+1},\xi_{k+1})$.
\STATE \textbf{If} $\rho_k\ne0$, evaluate $\nabla f(x_k,\xi_{k+1})$ and set
$g_{k+1}\gets \nabla f(x_{k+1},\xi_{k+1})+\rho_k\bigl(g_k-\nabla f(x_k,\xi_{k+1})\bigr)$.
\STATE $p_{k+1}\gets\frac12\left[(1-\tilde\alpha_k)p_k+
\sqrt{(1-\tilde\alpha_k)^2p_k^2+4\tilde\alpha_k\gamma_k\|g_{k+1}\|^2}\right]$.
\STATE $y_{k+\frac12}\gets y_k-\tilde\alpha_kp_{k+1}^{-1}g_{k+1}$, \quad $y_{k+1}\gets y_{k+\frac12}/\max\{1,2\|y_{k+\frac12}\|/R\}$.
\STATE $x_{k+1}^+\gets x_{k+1}-\eta_{k+1}p_{k+1}^{-1}g_{k+1}$.
\ENDFOR
\RETURN $x_K^+$.
\end{algorithmic}
\end{algorithm}

We retain the preconditioning structure of {\it Adam-SHANG-s}
\citep{yu2026adamshangconvergentadamtypemethod}, specialized here
to scalar preconditioning. The recursive estimate \eqref{eq:algorithmvrshang-g} replaces the fresh sampled gradient in both the search direction and the scalar preconditioner $p_kI$. The suffix \emph{-s} in the baseline's name denotes use of the current gradient in both the direction and the preconditioner.

Algorithm~\ref{alg:vr-shang} uses the polynomial steps and initialization
in Corollary~\ref{cor:polynomial-rate1}, with the parameter coupling
\eqref{eq:theorem-parameters1}.
Gradient feedback adapts the preconditioner $p_{k+1}$,
which scales both the auxiliary and correction updates.

These guarantees require $\mathcal B_R$ to contain a minimizer.
A sufficiently large $R$ ensures containment, but reduces the
feedback coefficient $\gamma_k=\tilde\alpha_k/R^2$.
Thus, $R$ controls both the projection domain and the strength
of gradient feedback. In practice, it may be treated as a separate
tuning parameter; the guarantees remain applicable provided
the required containment holds.

Both queries are formed before drawing the fresh sample. The first iteration uses one stochastic gradient evaluation; each later iteration uses two with the same fresh sample or minibatch. The total is $2K-1$ gradient-oracle calls.

\subsection{Gradient-estimation error recursion}
The energy argument follows {\it Adam-SHANG} \citep{yu2026adamshangconvergentadamtypemethod}. 
The additional task is to close a coupled bound: query motion controls the recursive estimation error, while the estimate drives the preconditioner that controls query motion.

Define the gradient-estimation error by
\begin{equation}\label{eq:paired-noise1}
U_k=g_k-\nabla f(x_k),\qquad k\ge1.
\end{equation}
Subtracting $\nabla f(x_{k+1})$ from \eqref{eq:algorithmvrshang-g} gives
\begin{equation}\label{eq:error-decomposition1}
U_{k+1}=\rho_kU_k+(1-\rho_k)Z_{k+1}+\rho_k(Z_{k+1}-Z_k^+),
\end{equation}
where the two noises use the same fresh sample:
\[
\begin{aligned}
Z_{k+1}=\nabla f(x_{k+1},\xi_{k+1})-\nabla f(x_{k+1}), \quad
Z_k^+=\nabla f(x_k,\xi_{k+1})-\nabla f(x_k).
\end{aligned}
\]
The term $\rho_kU_k$ propagates past error, while
$(1-\rho_k)Z_{k+1}$ damps fresh noise. The remaining term captures
noise changes between queries, controlled by mean-square
smoothness. The next lemma quantifies these contributions
through conditional moment bounds.

\begin{lemma}[Conditional error moments]\label{lem:recursion1}
For $k\ge1$, suppose that Assumptions~\ref{ass:unbiased-oracle}--\ref{ass:smoothness} hold, that $x_k,x_{k+1},g_k$ are $\F_k$-measurable, and that $\rho_k\in[0,1)$ is $\F_k$-measurable. Let $g_{k+1}$ follow \eqref{eq:algorithmvrshang-g}, with the error defined in \eqref{eq:paired-noise1}. Then
\begin{equation}\label{eq:error-mean1}
\E_kU_{k+1}=\rho_kU_k,
\end{equation}
and the conditional total variance satisfies
\begin{equation}\label{eq:split-innovation1}
\begin{aligned}
\operatorname{Var}_k(U_{k+1})
\le\bigl((1-\rho_k)\sigma+\rho_kL_\delta\|x_{k+1}-x_k\|\bigr)^2.
\end{aligned}
\end{equation}
\end{lemma}
\begin{proof}
The fresh sample is independent of $\F_k$. Unbiasedness therefore gives $\E_k(Z_{k+1}-\rho_kZ_k^+)=0$, so \eqref{eq:error-decomposition1} yields \eqref{eq:error-mean1}. Subtracting this conditional mean gives
$$
\begin{aligned}
\E_k\|U_{k+1}-\rho_kU_k\|^2
=\E_k\|(1-\rho_k)Z_{k+1}+\rho_k(Z_{k+1}-Z_k^+)\|^2.
\end{aligned}
$$
Since $x_k,x_{k+1}$ are $\F_k$-measurable and $\xi_{k+1}$
is independent of $\F_k$, Assumption~\ref{ass:smoothness} gives
\begin{equation}\label{eq:noise-smoothness}
\begin{aligned}
\E_k\|Z_{k+1}-Z_k^+\|^2
\le L_\delta^2\|x_{k+1}-x_k\|^2.
\end{aligned}
\end{equation}
The triangle inequality in conditional $L^2$ and the variance bound give
$$
\begin{aligned}
\sqrt{\E_k\|(1-\rho_k)Z_{k+1}+\rho_k(Z_{k+1}-Z_k^+)\|^2}
&\le |1-\rho_k|\sqrt{\E_k\|Z_{k+1}\|^2}
+|\rho_k|\sqrt{\E_k\|Z_{k+1}-Z_k^+\|^2}\\
&\le |1-\rho_k|\sigma+|\rho_k|L_\delta\|x_{k+1}-x_k\|.
\end{aligned}
$$
Squaring this bound proves \eqref{eq:split-innovation1}. 
\end{proof}

The conditional mean-square error exhibits a bias--variance tradeoff:
\begin{equation}\label{eq:tradeoff}
\begin{aligned}
\E_k\|U_{k+1}\|^2
&=\underbrace{\rho_k^2\|U_k\|^2}_{\text{squared conditional bias}}
+\underbrace{\E_k\|U_{k+1}-\rho_kU_k\|^2}_{\operatorname{Var}_k(U_{k+1})}\\
&\le \rho_k^2\|U_k\|^2 + \bigl((1-\rho_k)\sigma+\rho_kL_\delta\|x_{k+1}-x_k\|\bigr)^2.
\end{aligned}
\end{equation}
Increasing the memory weight retains more past error and damps the fresh-noise term,
while the noise-difference term also depends on query motion.
With $1-\rho_k=\tilde\alpha_k$, the query-motion bound
$\|x_{k+1}-x_k\|\le5R\tilde\alpha_k$ in
Lemma~\ref{lem:query-motion1} makes the variance contribution
of order $\tilde\alpha_k^2$, while the squared conditional bias
contracts by $(1-\tilde\alpha_k)^2$.

\subsection{Energy dissipation and self-normalization}

For $z=(x,y)$ and $P\succ0$, define
\begin{equation}\label{eq:energy1}
\cE(z,P)=f(x)-f(x^\star)+\frac12\|y-x^\star\|_P^2.
\end{equation}
The energy is the objective gap plus a metric-weighted distance to the minimizer. For the scalar version, set $P_k=p_kI$, $z_k^+=(x_k^+,y_k)$, and $z_k=(x_k,y_k)$. Split its change into a gradient-correction term and an interpolation--metric term:
\begin{equation}\label{eq:three-block-decomposition1}
\cE(z_{k+1}^+,P_{k+1})-\cE(z_k^+,P_k)
=\mathbb I_1+\mathbb I_2,
\end{equation}
where $\mathbb I_1 =\cE(z_{k+1}^+,P_{k+1})-\cE(z_{k+1},P_{k+1})$ and $\mathbb I_2=\cE(z_{k+1},P_{k+1})-\cE(z_k^+,P_k)$. We estimate the gradient-correction and interpolation–metric terms separately.

\begin{lemma}[Descent block]\label{lem:descent-block1}
Let $f$ be $L$-smooth and suppose that the correction step \eqref{eq:algorithmvrshang-xplus} satisfies
\begin{equation}\label{eq:descent-ratio1}
0<\eta_{k+1}\le\frac{p_{k+1}}L.
\end{equation}
Then, pathwise,
\begin{equation}\label{eq:descent-block-bound2}
\mathbb I_1\le-\frac{\eta_{k+1}p_{k+1}^{-1}}{2}\|g_{k+1}\|^2
+\eta_{k+1}p_{k+1}^{-1}\langle U_{k+1},g_{k+1}\rangle.
\end{equation}
\end{lemma}
\begin{proof}
The correction changes only $x$. By $L$-smoothness and $g_{k+1}=\nabla f(x_{k+1})+U_{k+1}$,
$$
\begin{aligned}
\mathbb I_1
&\le-\eta_{k+1}p_{k+1}^{-1}\langle\nabla f(x_{k+1}),g_{k+1}\rangle+\frac{L\eta_{k+1}^2p_{k+1}^{-2}}{2}\|g_{k+1}\|^2\\
&=-\eta_{k+1}p_{k+1}^{-1}\left(1-\frac{L\eta_{k+1}p_{k+1}^{-1}}{2}\right)\|g_{k+1}\|^2+\eta_{k+1}p_{k+1}^{-1}\langle U_{k+1},g_{k+1}\rangle.
\end{aligned}
$$
The step bound \eqref{eq:descent-ratio1} gives \eqref{eq:descent-block-bound2}.
\end{proof}

\begin{lemma}[Projected interpolation block]\label{lem:bregman-block1}
Suppose that $f$ is convex and that the updates \eqref{eq:algorithmvrshang} use $\gamma_k=\tilde\alpha_k/R^2$, with $y_k,x^\star\in\mathcal B_R$. Then, pathwise,
\begin{equation}\label{eq:bregman-block-bound1}
\begin{aligned}
\mathbb I_2\le{}&
-\tilde\alpha_k\cE(z_k^+,P_k)
-\tilde\alpha_k\langle U_{k+1},y_k-x^\star\rangle +\tilde\alpha_k^2p_{k+1}^{-1}
\|g_{k+1}\|^2.
\end{aligned}
\end{equation}
\end{lemma}
\begin{proof}
Projection does not increase the distance to $x^\star$. Using the unprojected step, the preconditioner update, and $\|y_k-x^\star\|\le R$ gives
$$
\begin{aligned}
&\frac{p_{k+1}}2\|y_{k+1}-x^\star\|^2-\frac{p_k}2\|y_k-x^\star\|^2 \\
&\le \frac{p_{k+1}}2\|y_{k+\frac{1}{2}}-x^\star\|^2-\frac{p_k}2\|y_k-x^\star\|^2 \\
&=-\frac{\tilde\alpha_kp_k}{2}\|y_k-x^\star\|^2-\tilde\alpha_k\langle g_{k+1},y_k-x^\star\rangle\\
&\quad+\frac{\tilde\alpha_k^2p_{k+1}^{-1}}{2}\|g_{k+1}\|^2
+\frac{\tilde\alpha_k\gamma_kp_{k+1}^{-1}}{2}\|y_k-x^\star\|^2\|g_{k+1}\|^2\\
&\le-\frac{\tilde\alpha_kp_k}{2}\|y_k-x^\star\|^2-\tilde\alpha_k\langle g_{k+1},y_k-x^\star\rangle+\tilde\alpha_k^2p_{k+1}^{-1}\|g_{k+1}\|^2.
\end{aligned}
$$
The second inequality uses $\gamma_k=\tilde\alpha_k/R^2$. Interpolation gives
$$\tilde\alpha_k(y_k-x^\star)=\tilde\alpha_k(x_{k+1}-x^\star)+(1-\tilde\alpha_k)(x_{k+1}-x_k^+).$$
Applying convexity to the two terms on the right yields
$$
\begin{aligned}
f(x_{k+1})-f(x_k^+)-\tilde\alpha_k\langle\nabla f(x_{k+1}),y_k-x^\star\rangle \le-\tilde\alpha_k\bigl(f(x_k^+)-f(x^\star)\bigr).
\end{aligned}
$$
Substituting $g_{k+1}=\nabla f(x_{k+1})+U_{k+1}$ into the weighted-distance estimate and adding the objective bound proves \eqref{eq:bregman-block-bound1}.
\end{proof}

The choice \(\eta_{k+1}=2\tilde\alpha_k^2\) cancels the two gradient-square terms in Lemmas~\ref{lem:descent-block1} and~\ref{lem:bregman-block1}. Substituting their bounds into \eqref{eq:three-block-decomposition1} gives
\begin{equation}\label{eq:basic-energy1}
	\cE(z_{k+1}^+,P_{k+1})-(1-\tilde\alpha_k)\cE(z_k^+,P_k)
	\le
	\langle U_{k+1},
	2\tilde\alpha_k^2p_{k+1}^{-1}g_{k+1}
	-\tilde\alpha_k(y_k-x^\star)\rangle.
\end{equation}
The pathwise bound \eqref{eq:basic-energy1} applies to both recursive
and nonrecursive estimates. 

The next lemma controls the remaining error pairing through pathwise
bounds on correction lengths and query motion. 
Here \emph{self-normalization} means that the same gradient estimate
sets both the direction and its scale: the feedback term ensures
$p_{k+1}\ge\tilde\alpha_k\|g_{k+1}\|/R$. Thus the normalized
correction remains bounded even when $\|g_{k+1}\|$ is large. 
Section~\ref{sec:mechanisms} discusses the roles of feedback and
recursive estimation.

\begin{lemma}[Self-normalization and energy bound]\label{lem:query-motion1}
Consider \eqref{eq:algorithmvrshang} with the correction coupling \(\eta_{k+1}=2\tilde\alpha_k^2\), $\gamma_k=\tilde\alpha_k/R^2$, $p_0>0$, and $x_0,x_0^+,y_0,x^\star\in\mathcal B_R$. Suppose that $0<\tilde\alpha_k<1$ for $k\ge0$ and
\begin{equation}\label{eq:alphacondition}
	(1-\tilde\alpha_k)\tilde\alpha_{k-1}\le\tilde\alpha_k, \qquad k \ge 1.
\end{equation}
Then, pathwise,
\begin{align}
\eta_{k+1}\|p_{k+1}^{-1}g_{k+1}\|&\le2R\tilde\alpha_k,\qquad k\ge0,
\label{eq:correction-length1}\\
\max\{\|x_k^+-x^\star\|,\|x_k-x^\star\|\}&\le3R,\qquad k\ge0,
\label{eq:uniform-iterates1}\\
\|x_{k+1}-x_k\|&\le5R\tilde\alpha_k,\qquad k\ge1.
\label{eq:query-motion1}
\end{align}
If $f$ is also convex and $L$-smooth and the correction step satisfies \eqref{eq:descent-ratio1}, then
\begin{equation}\label{eq:direct-energy1}
\begin{aligned}
\cE(z_{k+1}^+,P_{k+1})- (1-\tilde\alpha_k) \cE(z_k^+,P_k)
\le{}&3R\tilde\alpha_k\|U_{k+1}\|.
\end{aligned}
\end{equation}
If Assumptions~\ref{ass:unbiased-oracle}--\ref{ass:smoothness}
hold, the steps $\tilde\alpha_k$ are $\F_k$-measurable, and
$\rho_k=1-\tilde\alpha_k$ for $k\ge1$, then
\begin{equation}\label{eq:Ubound}
\E_k\|U_{k+1}\|^2
\le(1-\tilde\alpha_k)^2\|U_k\|^2
+\tilde\alpha_k^2(\sigma+5L_\delta R)^2,
\qquad k\ge1.
\end{equation}
\end{lemma}
\begin{proof}
The feedback term is nonnegative, so
\begin{equation}\label{eq:metric-decay-bound1}
p_{k+1}\ge(1-\tilde\alpha_k)p_k>0.
\end{equation}
Multiplying the preconditioner update by $p_{k+1}$ yields
$$
p_{k+1}^2
=(1-\tilde\alpha_k)p_kp_{k+1}
+\tilde\alpha_k^2R^{-2}\|g_{k+1}\|^2
\ge\tilde\alpha_k^2R^{-2}\|g_{k+1}\|^2.
$$
Taking square roots and using \(\eta_{k+1}=2\tilde\alpha_k^2\) proves the correction bound \eqref{eq:correction-length1}.

Projection keeps $y_k,x^\star\in\mathcal B_R$, so
$
\|y_k-x^\star\|\le\|y_k\|+\|x^\star\|\le R.
$
To bound the unprojected iterates, combine interpolation and correction:
$$
x_{k+1}^+=(1-\tilde\alpha_k)x_k^+
+\tilde\alpha_k y_k-\eta_{k+1}p_{k+1}^{-1}g_{k+1}.
$$
By \eqref{eq:correction-length1},
$$
\begin{aligned}
\|x_{k+1}^+\|
\le(1-\tilde\alpha_k)\|x_k^+\|
+\tilde\alpha_k\|y_k\|+\eta_{k+1}\|p_{k+1}^{-1}g_{k+1}\|
\le(1-\tilde\alpha_k)\|x_k^+\|+\tfrac{5R}{2}\tilde\alpha_k.
\end{aligned}
$$
Since $\|x_0^+\|\le R/2$, induction gives $\|x_k^+\|\le5R/2$. Interpolation gives the same bound for $\|x_k\|$. Together with $\|x^\star\|\le R/2$, these bounds prove \eqref{eq:uniform-iterates1} and give
$$
\|y_k-x_k\|\le\|y_k\|+\|x_k\|
\le3R.
$$
For $k\ge1$, the preceding correction implies
$$
x_{k+1}-x_k
=\tilde\alpha_k(y_k-x_k)
-(1-\tilde\alpha_k)\eta_kp_k^{-1}g_k.
$$
The correction bound and \eqref{eq:alphacondition} now give
$$
\begin{aligned}
\|x_{k+1}-x_k\|
\le\tilde\alpha_k\|y_k-x_k\|
+(1-\tilde\alpha_k)\eta_k\|p_k^{-1}g_k\|
\le3R\tilde\alpha_k
+2R(1-\tilde\alpha_k)\tilde\alpha_{k-1}
\le5R\tilde\alpha_k.
\end{aligned}
$$
Finally, \eqref{eq:correction-length1} and $\|y_k-x^\star\|\le R$ bound the error term in \eqref{eq:basic-energy1}:
\begin{equation}\label{eq:linear-error-bound1}
	\begin{aligned}
		\bigl\langle U_{k+1},
		\eta_{k+1}p_{k+1}^{-1}g_{k+1}-\tilde\alpha_k(y_k-x^\star)\bigr\rangle
		&\le\|U_{k+1}\|
		\bigl(\eta_{k+1}\|p_{k+1}^{-1}g_{k+1}\|
		+\tilde\alpha_k\|y_k-x^\star\|\bigr)\\
		&\le3R\tilde\alpha_k\|U_{k+1}\|.
	\end{aligned}
\end{equation}
Under the stated convexity, smoothness, and step conditions, \eqref{eq:basic-energy1} applies. Substituting \eqref{eq:linear-error-bound1} proves \eqref{eq:direct-energy1}. The conditional mean-square bound follows from
Lemma~\ref{lem:recursion1} and \eqref{eq:query-motion1}, using $1-\rho_k=\tilde\alpha_k$ and $\rho_k\le1$.
\end{proof}

\subsection{Convergence bounds and rates}
The preceding lemmas control one step; convergence requires controlling
how estimation errors accumulate over all iterations. We use weights matching
the contraction factors $1-\tilde\alpha_k$ and
$(1-\tilde\alpha_k)^2$ in the energy and mean-square error bounds:
\begin{equation}\label{eq:accelerated-weights1}
	A_0=1,\qquad
	A_{k+1} = A_k/(1-\tilde\alpha_k)=(1+\alpha_k)A_k.
\end{equation}
Weighting \eqref{eq:direct-energy1} and \eqref{eq:Ubound} by
$A_{k+1}$ and $A_{k+1}^2$, respectively, and summing separates
the initial-energy contribution from accumulated noise.
The following theorem gives a bound for feasible steps. The corollaries
then specialize it to polynomial and spectral steps.

\begin{theorem}[Stochastic energy bound]\label{thm:stochastic-objective1}
Suppose that $f$ satisfies \eqref{ass:convexity}, that Assumptions~\ref{ass:unbiased-oracle}--\ref{ass:smoothness} hold, and that $x^\star\in\mathcal B_R$. Initialize $x_0,x_0^+,y_0\in\mathcal B_R$, $\rho_0=0$, and $p_0>0$. Use the update rule \eqref{eq:algorithmvrshang} with $\F_k$-measurable steps $\tilde\alpha_k$ satisfying the following \emph{feasibility conditions} almost surely:
\begin{equation}\label{eq:step-feasibility1}
\begin{gathered}
0<\tilde\alpha_k\le
 \textstyle \min\left\{1/2,\sqrt{p_{k+1}/(2L)}\right\} \quad (k \ge 0), \quad \text{and} \quad 
(1-\tilde\alpha_k)\tilde\alpha_{k-1}\le\tilde\alpha_k \quad (k \ge 1).
\end{gathered}
\end{equation}
Choose the remaining parameters as
\begin{equation}\label{eq:theorem-parameters1}
\begin{gathered}
\alpha_k=\tilde\alpha_k/(1-\tilde\alpha_k),\qquad
\gamma_k=\tilde\alpha_k/R^2,\qquad
\eta_{k+1}=2\tilde\alpha_k^2,\qquad \rho_k=1-\tilde\alpha_k\quad(k\ge1).
\end{gathered}
\end{equation}
Then, for every $K\ge1$,
\begin{equation}\label{eq:stochastic-objective-bound1}
\begin{aligned}
\E\bigl[\cE(z_K^+,P_K)\bigr]
\le{}&\E[A_K^{-1}\cE(z_0^+,P_0)]\\
&+3R \textstyle\sum_{k=0}^{K-1}\sqrt{\E[\tilde\alpha_k^2A_K^{-2}]}\sqrt{4\sigma^2+(\sigma+5L_\delta R)^2
\textstyle\sum_{i=1}^{k}\E[A_{i+1}^2\tilde\alpha_i^2]}.
\end{aligned}
\end{equation}
The inner sum is zero for $k=0$.
\end{theorem}
\begin{proof}
The feasibility conditions and parameter coupling imply the correction-step bound \eqref{eq:descent-ratio1} and the hypotheses of Lemma~\ref{lem:query-motion1}, so \eqref{eq:direct-energy1} applies. Multiplying \eqref{eq:direct-energy1} by $A_{k+1}$ gives
$$
\begin{aligned}
A_{k+1}\cE(z_{k+1}^+,P_{k+1})-A_k\cE(z_k^+,P_k)
\le{}&3R\tilde{\alpha}_kA_{k+1}\|U_{k+1}\|.
\end{aligned}
$$
Summing over $k=0,\ldots,K-1$ and dividing by $A_K>0$ gives
\begin{equation}\label{eq:normalized-energy-bound1}
\begin{aligned}
\cE(z_K^+,P_K)
\le{}&A_K^{-1}\cE(z_0^+,P_0)
+3R \textstyle \sum_{k=0}^{K-1}\tilde{\alpha}_kA_{k+1}A_K^{-1}\|U_{k+1}\|.
\end{aligned}
\end{equation}
The initial sample gives $U_1=\nabla f(x_1,\xi_1)-\nabla f(x_1)$, hence $\E_0\|U_1\|^2\le\sigma^2$ by Assumption~\ref{ass:bounded-variance}. No bound on the initial query displacement $x_1-x_0$ is needed. 

Since $A_{k+1}$ is $\F_k$-measurable and
$A_{k+1}(1-\tilde\alpha_k)=A_k$, multiplying
\eqref{eq:Ubound} by $A_{k+1}^2$ yields
\begin{equation*}
	\E_k[A_{k+1}^2\|U_{k+1}\|^2]
	\le A_k^2\|U_k\|^2
	+A_{k+1}^2\tilde\alpha_k^2(\sigma+5L_\delta R)^2.
\end{equation*}
The initial mean-square gradient-estimation error bound and $A_1\le2$ imply $\E[A_1^2\|U_1\|^2]\le4\sigma^2$. Taking expectation and summing from $i=1$ to $k$ therefore gives
\begin{equation}\label{eq:weighted-error-moment1}
\E[A_{k+1}^2\|U_{k+1}\|^2]
\le4\sigma^2+(\sigma+5L_\delta R)^2
\textstyle \sum_{i=1}^{k}\E[A_{i+1}^2\tilde\alpha_i^2],\qquad k\ge0.
\end{equation}
Finally, Cauchy--Schwarz implies
$$
\begin{aligned}
\E\!\left[\tilde{\alpha}_kA_{k+1}A_K^{-1}\|U_{k+1}\|\right]
=\E\!\left[\tilde\alpha_kA_K^{-1}\|A_{k+1}U_{k+1}\|\right]\le\sqrt{\E[\tilde\alpha_k^2A_K^{-2}]}\,
\sqrt{\E[A_{k+1}^2\|U_{k+1}\|^2]}.
\end{aligned}
$$
No independence between $A_K$ and $U_{k+1}$ is needed. Taking expectation in \eqref{eq:normalized-energy-bound1} and applying \eqref{eq:weighted-error-moment1} proves \eqref{eq:stochastic-objective-bound1}. 

\end{proof}

The cap $1/2$ in \eqref{eq:step-feasibility1} is chosen for convenience:
it ensures $A_1\le2$, yielding the initial-error term $4\sigma^2$
in \eqref{eq:stochastic-objective-bound1}. Any fixed cap
$\bar\alpha\in[1/2,1)$ can be used instead, with this term replaced
by $\sigma^2/(1-\bar\alpha)^2$.
 
\begin{corollary}[Polynomial-step stochastic rate]\label{cor:polynomial-rate1}
Under the problem and initialization assumptions of
Theorem~\ref{thm:stochastic-objective1}, use the parameter coupling
\eqref{eq:theorem-parameters1} with
\begin{equation}\label{eq:polynomial-parameters1}
p_0=2L,\qquad \tilde\alpha_k=2/(k+4),\qquad k\ge0,
\end{equation}
Then the steps satisfy \eqref{eq:step-feasibility1}, and for every $K\ge1$,
\begin{equation}\label{eq:polynomial-tracking-rate1}
\E\|U_{K}\|^2
\le\frac4{K+3}(\sigma+\frac5{\sqrt2}L_\delta R)^2,
\end{equation}
\begin{equation}\label{eq:polynomial-rate1}
\E[\cE(z_K^+,P_K)]
\le\frac{6\E[\cE(z_0^+,P_0)]}{(K+2)(K+3)}
+\frac{8R\sigma+20\sqrt2\,L_\delta R^2}{\sqrt{K+3}}.
\end{equation}
Query-independent additive noise gives $L_\delta=0$; exact gradients
recover $O(1/K^2)$ convergence under \eqref{eq:polynomial-parameters1}. 
\end{corollary}
\begin{proof}
For the polynomial choice \eqref{eq:polynomial-parameters1}, the
weight recursion and \eqref{eq:metric-decay-bound1} give
\[
A_k=(k+2)(k+3)/6,\qquad p_k\ge 2LA_k^{-1}.
\]
Thus $p_{k+1}\ge2LA_{k+1}^{-1}\ge3L\tilde\alpha_k^2$.
Together with $0<\tilde\alpha_k\le1/2$ and
\[
(1-\tilde\alpha_k)\tilde\alpha_{k-1}
=\frac{2(k+2)}{(k+4)(k+3)}
\le\tilde\alpha_k,\qquad k\ge1,
\]
this verifies \eqref{eq:step-feasibility1}. Since $(\sigma+5L_\delta R)^2\le
2(\sigma+5L_\delta R/\sqrt2)^2$, induction in \eqref{eq:Ubound}
from $\E\|U_1\|^2\le\sigma^2$ gives, for $k\ge1$,
\[
\E\|U_{k+1}\|^2
\le\frac{4(\sigma+5L_\delta R/\sqrt2)^2}{(k+4)^2}
\left(\frac{(k+2)^2}{k+3}+2\right)
\le\frac{4(\sigma+5L_\delta R/\sqrt2)^2}{k+4}.
\]
This proves \eqref{eq:polynomial-tracking-rate1}.
Taking expectations in \eqref{eq:normalized-energy-bound1},
using Cauchy--Schwarz and $A_{k+1}\tilde\alpha_k=(k+3)/3$, yields
\[
\E[\cE(z_K^+,P_K)]
\le \E[\cE(z_0^+,P_0)]A_K^{-1}
+2RA_K^{-1}(\sigma+5L_\delta R/\sqrt2)
\sum_{k=0}^{K-1}\frac{k+3}{\sqrt{k+4}}.
\]
Finally, $\sqrt{1+t}\le1+t/2$ gives
\[
\frac{k+3}{\sqrt{k+4}}
\le\frac23\left[(k+3)\sqrt{k+4}-(k+2)\sqrt{k+3}\right].
\]
Thus the sum is at most $2(K+2)\sqrt{K+3}/3$, which gives
\eqref{eq:polynomial-rate1}.
\end{proof}

\begin{corollary}[Spectral-step feasibility and initial-energy decay]\label{rem:spectral-step1}
Under the problem and initialization assumptions of
Theorem~\ref{thm:stochastic-objective1}, use the parameter coupling
\eqref{eq:theorem-parameters1}.
The spectral rule
\begin{equation}\label{eq:square-root-step1}
\tilde\alpha_k=\min\left\{\frac13,\sqrt{\frac{p_k}{3L}}\right\}
\end{equation}
is predictable and satisfies \eqref{eq:step-feasibility1}, so the
bound \eqref{eq:stochastic-objective-bound1} applies. Moreover,
\begin{equation}\label{eq:square-root-clock1}
A_K\ge\left(1+K\min\left\{\frac16,
\sqrt{\frac{p_0}{12L}}\right\}\right)^2,\qquad K\ge0.
\end{equation}
Thus the initial-energy contribution decays as $O(1/K^2)$.
With exact gradients, the full energy satisfies
$\cE(z_K^+,P_K)\le\cE(z_0^+,P_0)/A_K=O(1/K^2)$.
\end{corollary}
\begin{proof}
The step is predictable because $p_k$ is $\F_k$-measurable.
Indeed, $p_{k+1}\ge(1-\tilde\alpha_k)p_k$ gives
\begin{equation}\label{eq:square-root-descent1}
2L\tilde\alpha_k^2\le\frac23p_k
\le(1-\tilde\alpha_k)p_k\le p_{k+1}.
\end{equation}
The same decay bound at the preceding step gives
\begin{equation}\label{eq:square-root-adjacent1}
\tilde\alpha_k\ge\tilde\alpha_{k-1}\sqrt{1-\tilde\alpha_{k-1}}
\ge \tilde\alpha_{k-1}/(1+\tilde\alpha_{k-1}),\qquad k\ge1.
\end{equation}
The cap preserves the first inequality, and the second follows from
$(1-t)(1+t)^2=1+t(1-t-t^2)\ge1$ for $0\le t\le1/3$.
Rearranging gives the adjacent-step condition.

For the initial-energy term, $A_kp_k\ge p_0$ gives
\[
A_{k+1}\tilde\alpha_k^2\ge\min\{1/9,p_0/(3L)\}.
\]
Hence
\[
\sqrt{A_{k+1}}-\sqrt{A_k}
=\frac{\sqrt{A_{k+1}}\tilde\alpha_k}{1+\sqrt{1-\tilde\alpha_k}}
\ge\min\left\{\frac16,\sqrt{\frac{p_0}{12L}}\right\}.
\]
Summation gives \eqref{eq:square-root-clock1}. With exact gradients,
the initial refresh gives $U_1=0$, and the recursion gives $U_k=0$
for all $k\ge1$. The energy rate follows from
\eqref{eq:normalized-energy-bound1}.
\end{proof}
The spectral rule lets the step respond to changes in the preconditioner.
The explicit stochastic rate \eqref{eq:polynomial-rate1} concerns
the polynomial choice.

\subsection{Roles of feedback and recursive estimation}\label{sec:mechanisms}
The proposed updates provide a common framework for HNAG-type
methods through different choices of the gradient estimator
and preconditioner feedback. The table summarizes these connections
under the projection and parameter couplings considered here:

\begin{center}
\renewcommand{\arraystretch}{1.35}
\resizebox{0.85\linewidth}{!}{%
\begin{tabular}{lll}
\toprule
Gradient estimate & No feedback $\gamma_k=0$
& Feedback $\gamma_k=\tilde\alpha_k/R^2$\\
\midrule
Exact & HNAG \citep{ChenLuo2019optimization}
& {\it Adam-HNAG-s} \citep{yu2026adamhnag}\\
Nonrecursive ($\rho_k=0$)
& SHANG \citep{yu2026shangrobuststochasticacceleration}
& {\it Adam-SHANG-s}
\citep{yu2026adamshangconvergentadamtypemethod}\\
Recursive & Not considered here
& {\it VR-Adam-SHANG} (this work)\\
\bottomrule
\end{tabular}%
}
\end{center}

\begin{enumerate}
\item \emph{No noise and no feedback: HNAG.}
Use exact gradients, so $U_k=0$, and set $\gamma_k=0$.
Let $A_k$ denote the weights defined in
\eqref{eq:accelerated-weights1}. The preconditioner then follows
a prescribed decay:
\[
p_{k+1}=(1-\tilde\alpha_k)p_k,\qquad p_k= p_0A_k^{-1}.
\]
For example, choose
\[
\tilde\alpha_k=2/(k+4),\qquad
\alpha_k=2/(k+2),\qquad
\eta_{k+1}=2\tilde\alpha_k^2,\qquad p_0=2L.
\]
Then $A_k=(k+2)(k+3)/6$ and
$p_k=12L/((k+2)(k+3))$, which ensure
$\eta_{k+1}\le p_{k+1}/L$. 
With no estimation error, it therefore gives
\[
\cE(z_K^+,P_K)
\le \cE(z_0^+,P_0)A_K^{-1}=O(1/K^2).
\]
This is the HNAG acceleration mechanism
\citep{ChenLuo2019optimization}.

\item \emph{Feedback with exact gradients: Adam-HNAG-s.}
Keep $U_k=0$ and set $\gamma_k=\tilde\alpha_k/R^2$.
The metric update now contains gradient feedback:
\[
p_{k+1}=(1-\tilde\alpha_k)p_k
+ R^{-2}\tilde\alpha_k^2p_{k+1}^{-1}
\|\nabla f(x_{k+1})\|^2.
\]
For fixed $p_k$ and $\tilde\alpha_k$, a larger gradient norm gives
a larger $p_{k+1}$. The spectral alternative in
Corollary~\ref{rem:spectral-step1} also lets the interpolation step
respond to this feedback, until its cap is reached. The accumulated factor
$A_K=\prod_{k=0}^{K-1}(1+\alpha_k)$ responds to these adaptive steps,
and the energy bound remains
\[
\cE(z_K^+,P_K)\le \cE(z_0^+,P_0)A_K^{-1}.
\]
This preserves the $O(1/K^2)$ guarantee while allowing gradient-dependent steps. This is the synchronous Adam-HNAG structure \citep{yu2026adamhnag}.

\item \emph{Sampling noise without feedback: SHANG.}
Set $\gamma_k=0$ and use a fresh, nonrecursive sampled gradient,
so $\rho_k=0$ and $\E_kU_{k+1}=0$.
The metric $p_{k+1}=(1-\tilde\alpha_k)p_k$ is known before sampling.
Under the descent condition, the energy estimate becomes
\[
\begin{aligned}
\E_k[\cE(z_{k+1}^+,P_{k+1})]
\le{}&(1-\tilde\alpha_k)\cE(z_k^+,P_k) +\eta_{k+1}p_{k+1}^{-1}\E_k\|U_{k+1}\|^2.
\end{aligned}
\]
Setting $\rho_k=0$ in \eqref{eq:tradeoff} gives
$\E_k\|U_{k+1}\|^2\le\sigma^2$, with no decay implied by this bound.
With persistent noise, deterministic contraction alone therefore
does not yield an $O(1/K^2)$ stochastic objective bound. 
SHANG obtains accelerated guarantees under multiplicative noise,
using noise-dependent parameter choices to absorb the variance
contribution \citep{yu2026shangrobuststochasticacceleration}.
\item \emph{Sampling noise with feedback: Adam-SHANG.}
Add gradient feedback to the third case by setting
$\gamma_k=\tilde\alpha_k/R^2$, while retaining the fresh sampled
gradient and $\rho_k=0$. Since the current gradient determines
both the direction and the metric, this is the synchronous variant
{\it Adam-SHANG-s} \citep{yu2026adamshangconvergentadamtypemethod}.
Conditional unbiasedness still cancels the interpolation-error term,
but the correction term now contains a random preconditioner:
\[
\begin{aligned}
\E_k[\cE(z_{k+1}^+,P_{k+1})]
\le{}&(1-\tilde\alpha_k)\cE(z_k^+,P_k) +\eta_{k+1}\E_k\langle U_{k+1},p_{k+1}^{-1}g_{k+1}\rangle.
\end{aligned}
\]
Unlike the third case, $p_{k+1}^{-1}$ cannot be taken outside the
conditional expectation. The original {\it Adam-SHANG}
analysis \citep{yu2026adamshangconvergentadamtypemethod} uses the lagged metric $p_k$, which is known before sampling. This permits
taking $p_k^{-1}$ outside the conditional expectation, but leaves
a variance term that need not decay.

\end{enumerate}

{\it VR-Adam-SHANG} adds recursive estimation to the fourth case,
trading conditional bias for reduced variance. By
Lemma~\ref{lem:recursion1},
\[
\E_k\|U_{k+1}\|^2
=\rho_k^2\|U_k\|^2+\operatorname{Var}_k(U_{k+1}).
\]
With $\rho_k=1-\tilde\alpha_k$, the inherited error contracts,
while self-normalization bounds query motion and makes the variance
term of order $\tilde\alpha_k^2$; see \eqref{eq:Ubound}.
Polynomial steps then give $O(1/K)$ mean-square estimation error
and an $O(1/\sqrt K)$ stochastic contribution to the objective bound.

\section{Diagonal VR-Adam-SHANG}\label{sec:vr-adam-shang}

The diagonal version of {\it VR-Adam-SHANG} uses the same queries and recursive gradient estimates as the scalar version, with a diagonal preconditioner and box projection.

\subsection{Algorithm}
Choose $R>0$ such that $\|x^\star\|_\infty\le R/2$, set $\mathcal H_R=[-R/2,R/2]^d$, and initialize $x_0,x_0^+,y_0\in\mathcal H_R$, $\rho_0 =0$, and diagonal $P_0\succ0$. Keep the query and estimator updates \eqref{eq:algorithmvrshang-x}--\eqref{eq:algorithmvrshang-g}, including $g_1=\nabla f(x_1,\xi_1)$. Replace the scalar preconditioner update by
\begin{subequations}\label{eq:algorithm}
\begin{equation}\label{eq:algorithm-P}
\frac{P_{k+1}-P_k}{\tilde\alpha_k}
=-P_k+\gamma_kP_{k+1}^{-1}G_{k+1}^2.
\end{equation}
\end{subequations}
In the auxiliary step and correction, replace $p_{k+1}^{-1}g_{k+1}$ by $P_{k+1}^{-1}g_{k+1}$. Replace projection onto $\mathcal B_R$ by projection onto $\mathcal H_R$ in the $P_{k+1}$-norm; since $P_{k+1}$ is diagonal, this clips each coordinate to $[-R/2,R/2]$. Here $G_{k+1}^2=\diag(g_{k+1}.^2)$ contains the squared gradient components. 

\begin{algorithm}[H]
\caption{{\it VR-Adam-SHANG}: diagonal version (changes to Algorithm~\ref{alg:vr-shang})}
\label{alg:vr-adam-shang-s}
\small
Use Algorithm~\ref{alg:vr-shang} with the following changes:
\begin{enumerate}
\item Replace $\mathcal B_R$ by $\mathcal H_R$ and initialize $P_0=2LI$.
\item Update $P_{k+1}$ by \eqref{eq:explicit-metric}, use $P_{k+1}^{-1}g_{k+1}$ in the auxiliary step and correction, and clip $y_{k+\frac{1}{2}}$ to $\mathcal H_R$.
\end{enumerate}
\end{algorithm}

The preconditioner remains positive definite and satisfies $P_{k+1}\succeq(1-\tilde\alpha_k)P_k$, hence $A_{k+1}P_{k+1}\succeq A_kP_k$ for the weights in \eqref{eq:accelerated-weights1}. Its explicit solution is
\begin{equation}\label{eq:explicit-metric}
	P_{k+1}=\frac12\left[(1-\tilde\alpha_k)P_k+
	\big((1-\tilde\alpha_k)^2P_k^2+4\tilde\alpha_k\gamma_kG_{k+1}^2 \big)^{1/2}\right].
\end{equation}

\subsection{Convergence analysis}
The energy decomposition and gradient-estimation error recursion
are unchanged. We record only the changes needed for diagonal
scaling: a directional descent condition, coordinatewise bounds
from box projection, and the resulting dimension factors.
 
\begin{lemma}[Preconditioned descent block]\label{lem:descent-block}
Suppose that $f$ is $L$-smooth and that $x_{k+1}^+=x_{k+1}-\eta_{k+1}P_{k+1}^{-1}g_{k+1}$. If $g_{k+1}\ne0$ and
\begin{equation}\label{eq:descent-ratio}
0<\eta_{k+1}\le q_{k+1}/L,
\qquad
q_{k+1}=\|g_{k+1}\|_{P_{k+1}^{-1}}^2/\|g_{k+1}\|_{P_{k+1}^{-2}}^2,
\end{equation}
then, pathwise,
\begin{equation}\label{eq:descent-block-bound}
\mathbb I_1\le-\frac{\eta_{k+1}}2
\|g_{k+1}\|_{P_{k+1}^{-1}}^2
+\eta_{k+1}\langle U_{k+1},g_{k+1}\rangle_{P_{k+1}^{-1}}.
\end{equation}
If $g_{k+1}=0$, the bound holds for every finite $\eta_{k+1}>0$.
\end{lemma}
The proof follows the scalar argument in Lemma~\ref{lem:descent-block1}, replacing $p_{k+1}^{-1}g_{k+1}$ by $P_{k+1}^{-1}g_{k+1}$: the smoothness term is $L\eta_{k+1}^2\|g_{k+1}\|_{P_{k+1}^{-2}}^2/2$, which \eqref{eq:descent-ratio} bounds by $\eta_{k+1}\|g_{k+1}\|_{P_{k+1}^{-1}}^2/2$.

The ratio $q_{k+1}$ is a weighted average of the diagonal entries of $P_{k+1}$, with weights proportional to $((P_{k+1}^{-1}g_{k+1})_i)^2$. Hence
\begin{equation}\label{eq:ratio-spectrum}
\lambda_{\min}(P_{k+1})\le q_{k+1}\le\lambda_{\max}(P_{k+1}).
\end{equation}
The descent condition therefore depends on both the preconditioner
and the update direction.

For the interpolation block, Lemma~\ref{lem:bregman-block1}
applies coordinatewise: replace
$p_{k+1}^{-1}\|g_{k+1}\|^2$ by
$\|g_{k+1}\|_{P_{k+1}^{-1}}^2$ in
\eqref{eq:bregman-block-bound1}. Box projection gives
$|(y_k-x^\star)_i|\le R$, which bounds each coordinate of the
feedback term; the convexity argument is unchanged. Under
\eqref{eq:descent-ratio}, combining this estimate with
\eqref{eq:descent-block-bound} and choosing
$\eta_{k+1}=2\tilde\alpha_k^2$ gives
\begin{equation}\label{eq:basic-energy}
	\cE(z_{k+1}^+,P_{k+1})-(1-\tilde\alpha_k)\cE(z_k^+,P_k)
	\le
	\left\langle U_{k+1},
	2\tilde\alpha_k^2P_{k+1}^{-1}g_{k+1}
	-\tilde\alpha_k(y_k-x^\star)\right\rangle.
\end{equation}
To control the remaining error pairing, we bound the correction and query motion using a coordinatewise version of the scalar argument.

\begin{lemma}[Correction, iterate, and query bounds]\label{lem:query-motion}
Consider the diagonal updates above, with $\eta_{k+1}=2\tilde\alpha_k^2$, $\gamma_k=\tilde\alpha_k/R^2$, $P_0\succ0$, $x_0,x_0^+,y_0\in\mathcal H_R$, and $x^\star\in\mathcal H_R$. Suppose that $0<\tilde\alpha_k<1$ for $k\ge0$ and the adjacent-step condition \eqref{eq:alphacondition} holds. Then, pathwise,
\begin{align}
\eta_{k+1}\|P_{k+1}^{-1}g_{k+1}\|_\infty
&\le2R\tilde\alpha_k,\qquad k\ge0,\label{eq:correction-length}\\
\max\{\|x_k^+-x^\star\|_\infty,\|x_k-x^\star\|_\infty\}
&\le3R,\qquad k\ge0,\label{eq:uniform-iterates}\\
\|x_{k+1}-x_k\|&\le5\sqrt d\,R\tilde\alpha_k,\qquad k\ge1.
\label{eq:query-motion}
\end{align}
\end{lemma}
These bounds follow by applying the self-normalization and iterate bounds in Lemma~\ref{lem:query-motion1} coordinatewise, then using $\|v\|\le\sqrt d\,\|v\|_\infty$ for the query bound.

Each coordinate of the vector paired with $U_{k+1}$ in \eqref{eq:basic-energy} has magnitude at most $3R\tilde\alpha_k$. Thus the pairing is at most $3R\tilde\alpha_k\|U_{k+1}\|_1$. Using $\|U_{k+1}\|_1\le\sqrt d\,\|U_{k+1}\|$ gives
\begin{equation}\label{eq:linear-error-bound}
\left\langle U_{k+1},
\eta_{k+1}P_{k+1}^{-1}g_{k+1}
-\tilde\alpha_k(y_k-x^\star)\right\rangle
\le3\sqrt d\,R\tilde\alpha_k\|U_{k+1}\|.
\end{equation}

Combining \eqref{eq:basic-energy} with \eqref{eq:linear-error-bound} and summing with the weights $A_{k+1}$ from \eqref{eq:accelerated-weights1}, as in the proof of Theorem~\ref{thm:stochastic-objective1}, yields
\begin{equation}\label{eq:cumulative-energy}
	A_K\cE(z_K^+,P_K)
	\le\cE(z_0^+,P_0)+3\sqrt d\,R
	\textstyle \sum_{k=0}^{K-1}(A_{k+1}-A_k)\|U_{k+1}\|.
\end{equation}
Controlling the gradient-estimation error $U_{k+1}$ then gives the following convergence theorem.

\begin{theorem}[Stochastic energy bound]\label{thm:diagonal-energy}
Suppose that $f$ satisfies \eqref{ass:convexity}, that Assumptions~\ref{ass:unbiased-oracle}--\ref{ass:smoothness} hold, and that $x^\star\in\mathcal H_R$. Initialize $x_0,x_0^+,y_0\in\mathcal H_R$, $\rho_0=0$, and a deterministic diagonal matrix $P_0\succ0$. Use the diagonal updates with $\F_k$-measurable steps satisfying, almost surely,
\begin{equation}\label{eq:diagonal-feasibility}
\begin{gathered}
0<\tilde\alpha_k\le \textstyle \min\left\{1/2,
\sqrt{q_{k+1}/(2L)}\right\}\quad(k\ge0),\qquad
(1-\tilde\alpha_k)\tilde\alpha_{k-1}\le\tilde\alpha_k\quad(k\ge1),
\end{gathered}
\end{equation}
where $q_{k+1}$ is the directional ratio in \eqref{eq:descent-ratio}.
The restriction involving $q_{k+1}$ is required only when $g_{k+1}\ne0$. Set
\begin{equation}\label{eq:theorem-coupling}
\alpha_k=\tilde\alpha_k/(1-\tilde\alpha_k),\qquad
\gamma_k=\tilde\alpha_k/R^2,\qquad
\eta_{k+1}=2\tilde\alpha_k^2,\qquad
\rho_k=1-\tilde\alpha_k\quad(k\ge1).
\end{equation}
Define $A_k$ by \eqref{eq:accelerated-weights1}. Then, for every $K\ge1$,
\begin{equation}\label{eq:stochastic-diagonal-bound}
\begin{aligned}
\E[\cE(z_K^+,P_K)]
\le{}&\E[A_K^{-1}\cE(z_0^+,P_0)]\\
&+3\sqrt d\,R \textstyle \sum_{k=0}^{K-1}
\sqrt{\E[\tilde\alpha_k^2A_K^{-2}]}\,
\sqrt{4\sigma^2+(\sigma+5\sqrt d\,L_\delta R)^2
\textstyle \sum_{i=1}^k\E[A_{i+1}^2\tilde\alpha_i^2]}.
\end{aligned}
\end{equation} 
The inner sum is zero for $k=0$.
In particular, $P_0=2LI$ and $\tilde\alpha_k=2/(k+4)$ satisfy
\eqref{eq:diagonal-feasibility}. The mean-square gradient-estimation
error satisfies \eqref{eq:polynomial-tracking-rate1} with $R$
replaced by $\sqrt d\,R$, and, for every $K\ge1$,
\begin{equation}\label{eq:polynomial-rate}
\E[\cE(z_K^+,P_K)]
\le\frac{6\E[\cE(z_0^+,P_0)]}{(K+2)(K+3)}
+\frac{8\sqrt d\,R\sigma+20\sqrt2\,dL_\delta R^2}{\sqrt{K+3}}.
\end{equation}
\end{theorem} 
\begin{proof}
The error recursion in Lemma~\ref{lem:recursion1} is unchanged. Substituting the query bound \eqref{eq:query-motion} gives
$$
\E_k\|U_{k+1}\|^2\le(1-\tilde\alpha_k)^2\|U_k\|^2
+\tilde\alpha_k^2(\sigma+5\sqrt d\,L_\delta R)^2,\qquad k\ge1.
$$
Together with $\E_0\|U_1\|^2\le\sigma^2$ and \eqref{eq:cumulative-energy}, this is the scalar proof with $R$ replaced by $\sqrt d\,R$.

For the polynomial choice, feasibility follows as in
Theorem~\ref{thm:stochastic-objective1}, using
$\lambda_{\min}(P_{k+1})\ge2L/A_{k+1}\ge3L\tilde\alpha_k^2$
and \eqref{eq:ratio-spectrum}. The scalar polynomial-step argument
with $R$ replaced by $\sqrt d\,R$ then gives the explicit error
bound and \eqref{eq:polynomial-rate}.
\end{proof}

\begin{remark}[Spectral steps as an alternative]\label{rem:diagonal-spectral}
The general bound \eqref{eq:stochastic-diagonal-bound} also applies
to other predictable steps satisfying \eqref{eq:diagonal-feasibility}.
For any deterministic diagonal $P_0\succ0$, one such choice is
\begin{equation}\label{eq:diagonal-spectral-step}
\tilde\alpha_k=\min\left\{\frac13,
\sqrt{\frac{\lambda_{\min}(P_k)}{3L}}\right\}.
\end{equation}
Indeed, when $g_{k+1}\ne0$,
\[
q_{k+1}\ge\lambda_{\min}(P_{k+1})
\ge(1-\tilde\alpha_k)\lambda_{\min}(P_k)
\ge3L(1-\tilde\alpha_k)\tilde\alpha_k^2
\ge2L\tilde\alpha_k^2.
\]
The adjacent-step condition and the $O(1/K^2)$ initial-energy
bound follow from Corollary~\ref{rem:spectral-step1} with
$p_k$ replaced by $\lambda_{\min}(P_k)$.
\end{remark}

\subsection{VR-Adam-SHANG without recursive estimation}\label{sec:adam-shang-s}

The energy and motion bounds also hold when $\rho_k=0$, so $g_{k+1}=\nabla f(x_{k+1},\xi_{k+1})$. This specialization of {\it VR-Adam-SHANG} has the nonrecursive {\it Adam-SHANG-s} structure \citep{yu2026adamshangconvergentadamtypemethod}, with the box projection and polynomial steps used here. 

\begin{corollary}[Polynomial-step bound without recursive estimation]
\label{cor:synchronous-polynomial}
Under the assumptions of Theorem~\ref{thm:diagonal-energy}, use
$P_0=2LI$, $\tilde\alpha_k=2/(k+4)$, and the initialization and
parameter couplings stated there, except that $\rho_k=0$ for all
$k\ge0$. Then, for every $K\ge1$,
\begin{equation}\label{eq:synchronous-polynomial-rate}
\E[f(x_K^+)-f(x^\star)] 
\le\frac{6\E[\cE(z_0^+,P_0)]}{(K+2)(K+3)}
+\frac{8R\sqrt{6d}\sqrt{9L^2dR^2+\sigma^2}}
{3\sqrt{K+2}}.
\end{equation}
\end{corollary}
\begin{proof}
The feasibility and motion bounds for the polynomial choice in
Theorem~\ref{thm:diagonal-energy} do not depend on the memory schedule,
so \eqref{eq:basic-energy} still applies. Here
$U_{k+1}=g_{k+1}-\nabla f(x_{k+1})$ is fresh gradient noise with
$\E_kU_{k+1}=0$. Taking conditional expectations cancels the
pairing with the $\F_k$-measurable vector $y_k-x^\star$, giving
\[
\E_k[\cE(z_{k+1}^+,P_{k+1})]
-(1-\tilde\alpha_k)\cE(z_k^+,P_k)
\le2\tilde\alpha_k^2
\E_k\langle U_{k+1},P_{k+1}^{-1}g_{k+1}\rangle.
\]
It remains to control the correlation on the right, since the
preconditioner uses the current gradient.
Fix $k$ and a coordinate $i$, and condition on $\F_k$. Substituting $\gamma_k=\tilde\alpha_k/R^2$ into the explicit preconditioner update \eqref{eq:explicit-metric} defines the scalar map
$$
\psi_{k,i}(t)=
2t \left[(1-\tilde\alpha_k)(P_k)_{ii}
+\sqrt{(1-\tilde\alpha_k)^2(P_k)_{ii}^2
+4\tilde\alpha_k^2t^2/R^2} \right]^{-1},\qquad t\in\R.
$$
It maps a gradient component to its preconditioned value:
$
(P_{k+1}^{-1}g_{k+1})_i=\psi_{k,i}((g_{k+1})_i).
$
Here $P_k$ and $\tilde\alpha_k$ are fixed under conditioning. The even denominator makes $\psi_{k,i}$ odd, and differentiation shows that it is strictly increasing. 

Let $X=(g_{k+1})_i$. Since $|\psi_{k,i}(X)|=\psi_{k,i}(|X|)$, Chebyshev's integral inequality for the increasing functions $t$ and $\psi_{k,i}(t)$ on $[0,\infty)$ gives
$$
\begin{aligned}
\left|\E_kX\,\E_k[\psi_{k,i}(X)]\right|
\le\E_k|X|\,\E_k[\psi_{k,i}(|X|)]\le\E_k[|X|\psi_{k,i}(|X|)]
=\E_k[X\psi_{k,i}(X)].
\end{aligned}
$$
Unbiasedness gives $(U_{k+1})_i=X-\E_kX$ and $X\psi_{k,i}(X) = (P^{-1}_{k+1})_{ii} (g_{k+1})_i^2$, so
$$
\begin{aligned}
\E_k[(U_{k+1})_i(P_{k+1}^{-1}g_{k+1})_i]
=\E_k[X\psi_{k,i}(X)]-\E_kX\,\E_k[\psi_{k,i}(X)]
\le2\E_k (P^{-1}_{k+1})_{ii} (g_{k+1})_i^2.
\end{aligned}
$$
No symmetry assumption on the noise distribution is needed. Summing over coordinates and applying the conditional energy estimate gives
$$
\E_k[\cE(z_{k+1}^+,P_{k+1})]
-(1-\tilde\alpha_k)\cE(z_k^+,P_k)
\le4\tilde\alpha_k^2\E_k\|g_{k+1}\|_{P_{k+1}^{-1}}^2.
$$
The preconditioner update gives $A_{k+1}P_{k+1}-A_kP_k=A_{k+1}\tilde\alpha_k^2P_{k+1}^{-1}G_{k+1}^2/R^2$. Taking traces and summing the energy bound therefore yields
$$
A_K\E[\cE(z_K^+,P_K)]
\le\E[\cE(z_0^+,P_0)]
+4R^2\E[\Tr(A_KP_K-P_0)].
$$
Factoring the difference of squares and using $A_{k+1}-A_k=A_{k+1}\tilde\alpha_k$ in the preconditioner update identity gives, coordinatewise,
$$
\begin{aligned}
A_{k+1}^2(P_{k+1})_{ii}^2-A_k^2(P_k)_{ii}^2
& =(A_{k+1}-A_k)^2(g_{k+1})_i^2R^{-2} \left(1+\tfrac{A_k(P_k)_{ii}}{A_{k+1}(P_{k+1})_{ii}}\right)\\
& \le2R^{-2}(A_{k+1}-A_k)^2(g_{k+1})_i^2.
\end{aligned}
$$
The inequality uses $A_kP_k\preceq A_{k+1}P_{k+1}$, so the ratio in parentheses is at most one. Lemma~\ref{lem:query-motion} gives $\|x_{k+1}-x^\star\|^2\le d\|x_{k+1}-x^\star\|_\infty^2\le9dR^2$. Unbiasedness eliminates the gradient--noise cross term; $L$-smoothness and $\nabla f(x^\star)=0$ then give
$$
\begin{aligned}
\E_k\|g_{k+1}\|^2
=\|\nabla f(x_{k+1})\|^2+\E_k\|U_{k+1}\|^2
\le L^2\|x_{k+1}-x^\star\|^2+\sigma^2
\le9L^2dR^2+\sigma^2.
\end{aligned}
$$
Summing the coordinatewise preconditioner bound, using $\sqrt{u+v}-\sqrt u\le\sqrt v$, and applying Cauchy--Schwarz therefore gives
$$
\E[\Tr(A_KP_K-P_0)]
\le\sqrt{2d}R^{-1}\sqrt{9L^2dR^2+\sigma^2}
\sqrt{\textstyle \sum_{k=0}^{K-1}(A_{k+1}-A_k)^2}.
$$
Finally, the polynomial weights satisfy
$$
A_K=(K+2)(K+3)/6,\quad
A_{k+1}-A_k=(k+3)/3,\quad
A_K^{-2}\textstyle \sum_{k<K}(A_{k+1}-A_k)^2
\le\tfrac4{3(K+2)}.
$$
Substituting into the energy bound proves \eqref{eq:synchronous-polynomial-rate}.
\end{proof}

The nonrecursive bound in Corollary~\ref{cor:synchronous-polynomial}
has the same $O(1/\sqrt K)$ stochastic order as the recursive bound
in Theorem~\ref{thm:diagonal-energy}, but differs in its noise dependence.
It uses the full gradient second moment, introducing a
noise-independent contribution of order $LdR^2/\sqrt K$.
The corresponding term in the recursive bound instead depends on
$\sigma$ and $L_\delta$. Thus the distinction between these bounds
lies in their noise dependence rather than their convergence orders.

Scalar preconditioning changes only the length of $g_{k+1}$,
whereas diagonal preconditioning rescales each coordinate separately
and can change its direction. This allows different coordinates
to adapt to heterogeneous gradient magnitudes.
Both versions have the same convergence orders; their distinction
is directional adaptation rather than a different rate order.
The dimension factors in the diagonal bounds arise from converting
box bounds to Euclidean bounds. Since $\mathcal B_R$ is a ball
and $\mathcal H_R$ is a box, constants at the same numerical $R$
are not directly comparable.

\section{Numerical Experiments}\label{sec:experiments}

We test the scalar and diagonal versions of {\it VR-Adam-SHANG} (Algorithms~\ref{alg:vr-shang} and~\ref{alg:vr-adam-shang-s}).
The baselines are projected SHANG \citep{yu2026shangrobuststochasticacceleration}, STORM \citep{cutkosky2019momentum}, Adam \citep{kingma2015adam}, {\it Adam-SHANG-s} \citep{yu2026adamshangconvergentadamtypemethod}, and $\mu^2$-ExtraSGD \citep{dahan2025mu2sgd}. The comparison with nonrecursive {\it Adam-SHANG-s} isolates recursive estimation under the same polynomial schedule and diagonal geometry.

\medskip

\noindent \textbf{Initialization and parameters.}
All methods evaluate their first minibatch gradient at $\mathbf1$.
SHANG-based methods set $x_0=\mathbf1$ and $y_0=0$, with $x_0^+=x_0$
for their adaptive variants.
Both proposed versions follow the initialization and parameter
choices in Algorithms~\ref{alg:vr-shang} and~\ref{alg:vr-adam-shang-s}.
Nonrecursive {\it Adam-SHANG-s} follows
Corollary~\ref{cor:synchronous-polynomial}. Projected SHANG uses $\alpha_k=2/(k+1)$, $p_k=L\alpha_k^2$, and $\eta_k=1/L$.
The legends abbreviate the proposed versions as ``Ours (scalar)''
and ``Ours (diagonal).''
Projection radii are fixed before the runs and pass the posterior
checks in Appendix~\ref{app:radius-verification}.

Adam initializes both moment vectors to zero and uses moment parameters
$(0.9,0.999)$, bias correction, $\varepsilon=10^{-8}$, and learning rate
$\ell_0/\sqrt{k+1}$. 
STORM uses $\eta_k=\kappa/(w+(k+1)\sigma^2)^{1/3}$ and recursive
coefficient $1-c\eta_k^2$, with $\kappa,w$ determined by $L,\sigma^2,c$
as in~\citep[Appendix~B]{cutkosky2019momentum}.
For $\mu^2$-ExtraSGD, we use Algorithm~2 in~\citep{dahan2025mu2sgd},
with weights $\alpha_k=k+1$, recursive coefficients $\beta_k=1/\alpha_k$,
step $\eta=\theta/(2L)$, and weighted-average output. Its extragradient
base point starts at $\mathbf1$, the accumulated averaging weight starts
at zero, and its projection ball is the same as that of scalar
{\it VR-Adam-SHANG}.

 For each dataset and noise setting, we tune Adam's $\ell_0$, STORM's $c$, and the step factor $\theta$ of $\mu^2$-ExtraSGD by minimizing the mean terminal objective at the comparison budget over three independent pilot runs. Selected parameters appear in the legends. Appendix~\ref{app:common-implementation} gives the grids and seeds.

\medskip
\noindent \textbf{Gradient budgets and evaluation.}
We use a batch size of $B=64$ throughout and measure computational
cost in sample-gradient evaluations. For $T$ updates, the costs are
$BT$ for SHANG and Adam, $B(T+1)$ for {\it Adam-SHANG-s},
$B(2T-1)$ for STORM, $2BT$ for both proposed methods, and
$B(3T-1)$ for $\mu^2$-ExtraSGD. 

We compare terminal objective gaps under a common budget of $3BK$
sample-gradient evaluations, with $K=2048$. Each method performs as
many complete updates as this budget permits, and its objective gap
is evaluated at the final output. Gradient-estimation mean-square errors are compared separately over
the first $K$ iterations of each method.

All methods share sample streams. Reference solutions and
exact-gradient diagnostics are used only for evaluation and
excluded from the sample-gradient budget. All curves show means over ten runs with
pointwise $95\%$ bootstrap confidence intervals.

\subsection{Graph-signal denoising}\label{sec:graph-denoising}
We consider graph-Laplacian Tikhonov denoising \citep{shuman2013emerging} on a $12\times12$ grid:
\begin{equation}\label{eq:experiment-objective}
	f(x)=\frac12 \E \|x-z\|_D^2 +\frac12x^\top L_Gx.
\end{equation}
Here $x\in\R^{144}$ is the estimated signal, $L_G$ is the graph Laplacian, and $D\succ0$ is diagonal. Its entries vary by a factor of $1000$, with $\Tr(D)=144$. Observations have the form $z=b+\nu D^{-1/2}\zeta$, where $b$ is a fixed signal and $\zeta\sim N(0,I)$. The data-generation procedure described in
Appendix~\ref{app:graph-details} yields an objective with $L\approx11.072$. The noise is additive, so $L_\delta=0$ and same-sample gradient differences are noise-free.  We choose $R=2\|x_0\|= 24$ for the scalar methods and $R=2\|x_0\|_\infty=2$ for the diagonal methods. We test $\nu\in\{0.1,1\}$.

\begin{figure}[!htb]
	\centering
	\includegraphics[width=0.495\linewidth]{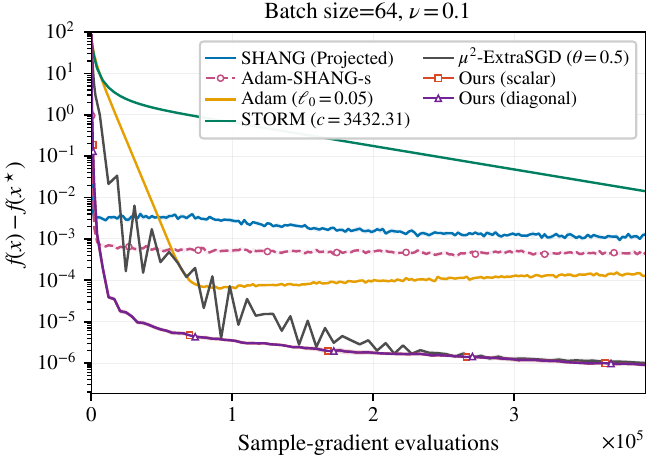}\hfill%
	\includegraphics[width=0.495\linewidth]{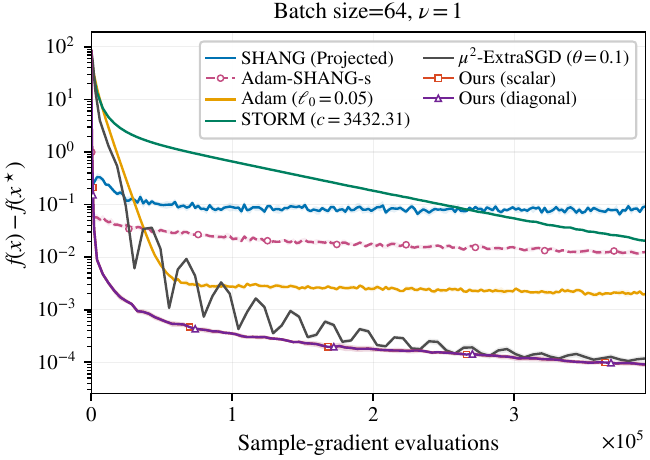}
	\caption{\small Objective gaps for graph-signal denoising at
		$\nu=0.1$ (left) and $\nu=1$ (right).
		The scalar and diagonal curves nearly overlap.}
	\label{fig:convex-comparison}
	\vspace{-0.8em}
\end{figure}

\begin{figure}[!htb]
	\centering
	\includegraphics[width=0.495\linewidth]{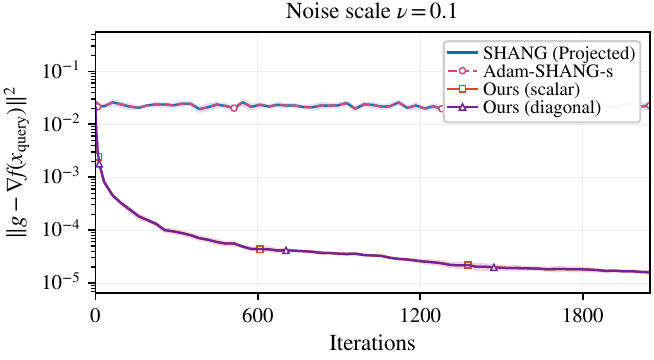}\hfill%
	\includegraphics[width=0.495\linewidth]{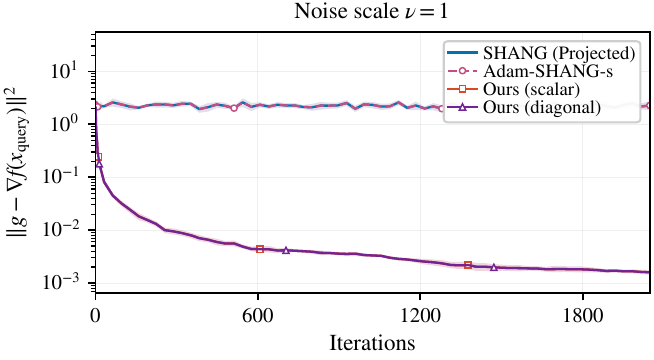}
	\caption{\small Mean-square gradient-estimation errors for graph-signal denoising. The SHANG and {\it Adam-SHANG-s} curves nearly overlap, as do the scalar and diagonal curves.}
	\label{fig:convex-tracking}
	\vspace{-0.8em}
\end{figure}

Both versions attain lower mean terminal gaps than the baselines at
both noise levels, with a modest advantage over $\mu^2$-ExtraSGD
and a larger advantage over nonrecursive {\it Adam-SHANG-s}
(Figure~\ref{fig:convex-comparison}). Mean-square estimation errors,
$\E\|g-\nabla f(x)\|^2$ at each method's gradient query, decrease
for recursive methods but remain at the minibatch-noise level for
nonrecursive methods (Figure~\ref{fig:convex-tracking}).

\subsection{Logistic regression on synthetic data}\label{sec:sparse-logistic}
We test query-dependent noise using sparse logistic regression on
$N=12000$ observations $(a_i,b_i)\in\R^{60}\times\{-1,1\}$:
\begin{equation}\label{eq:sparse-logistic-objective}
	f(x)=\tfrac1N\textstyle\sum_{i=1}^N
	\log\bigl(1+\exp(-b_i a_i^\top x)\bigr)
	+\frac{\lambda}{2}\|x\|^2,
	\qquad \lambda=10^{-4}.
\end{equation}
Feature frequencies (fractions of observations with a nonzero entry)
vary approximately geometrically, with a largest-to-smallest ratio
of $99.2$. The $72000$ nonzero entries independently take values
$\pm1$, without standardization. We use $L=4.25+\lambda$,
$R=2\sqrt{60}$ for scalar methods, and $R=2$ for diagonal methods;
Appendix~\ref{app:logistic-details} gives the data and oracle details.

\begin{figure}[!htb]
	\centering
	\includegraphics[width=0.495\linewidth]{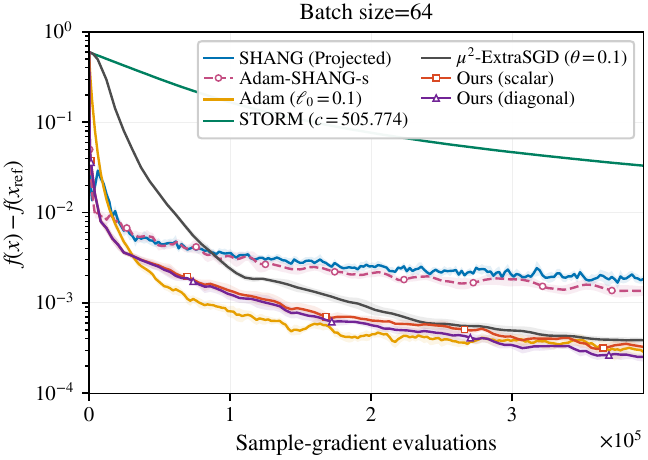}\hfill%
	\includegraphics[width=0.495\linewidth]{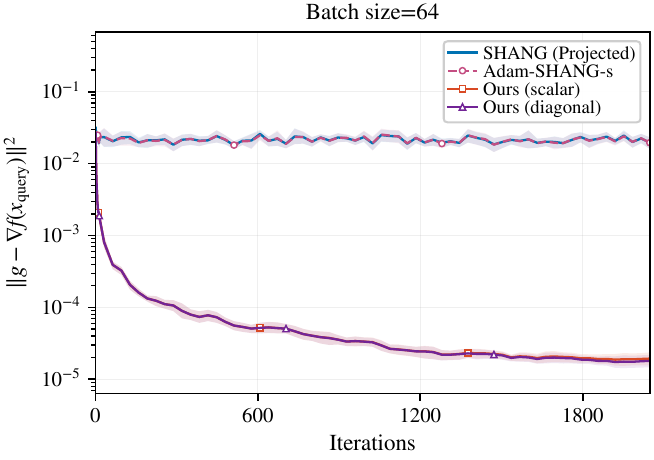}
	\caption{\small Objective gaps (left) and mean-square gradient-estimation
		errors (right) for logistic regression.}
	\label{fig:sparse-logistic}
	\vspace{-0.8em}
\end{figure}

Both versions initially reduce objective gaps faster than Adam and
continue improving after its gap levels off. The diagonal version
achieves the lowest mean terminal gap. Recursive estimation also
reduces gradient-estimation errors beyond the additive-noise setting
(Figure~\ref{fig:sparse-logistic}).

\begin{figure}[!htb]
	\centering
	\includegraphics[width=0.495\linewidth]{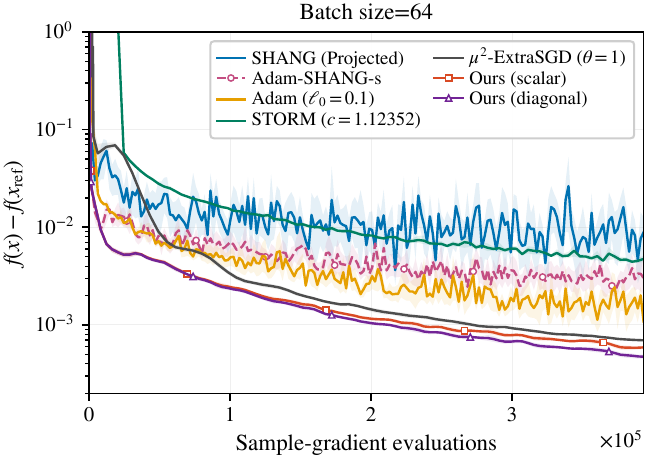}\hfill%
	\includegraphics[width=0.495\linewidth]{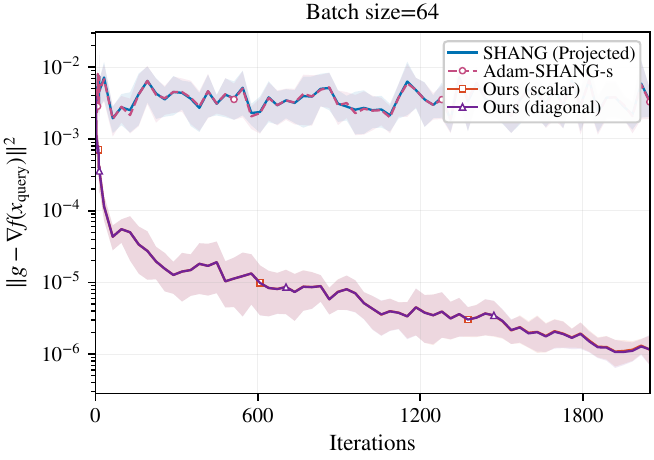}
	\caption{\small Objective gaps (left, clipped above $1$) and mean-square gradient-estimation errors (right) for CIFAR-10.}
	\label{fig:cifar10-logistic}
	\vspace{-1em}
\end{figure}

\subsection{Logistic regression on real data}\label{sec:cifar10-logistic}
We next use all $N=10000$ cat and dog training images from
CIFAR-10 \citep{Krizhevsky2009}, with $5000$ images per class.
Each RGB image is flattened into $d=3072$ raw pixel features and
normalized to unit Euclidean norm; labels are $-1$ for cats and
$1$ for dogs. We minimize \eqref{eq:sparse-logistic-objective}
with the same $\lambda=10^{-4}$, without an intercept,
augmentation, or learned features. The resulting problem is convex.

We use $R=2\sqrt{3072}\approx110.851$ for scalar methods and
$R=2$ for diagonal methods. A global smoothness bound is
$L\approx0.200314$, and STORM uses $\sigma^2=1/B$.

Figure~\ref{fig:cifar10-logistic} shows lower mean terminal gaps
for both versions than for all baselines, with the diagonal version
attaining the lowest gap. Their mean-square gradient-estimation
errors also decrease under recursion.

\vspace{-0.6em}
\section{Conclusion}\label{sec:conclusion}
We developed scalar and diagonal versions of {\it VR-Adam-SHANG} by combining
recursive gradient estimation with Hessian-driven acceleration and
scalar or coordinatewise Adam-style preconditioning. By refining the
gradient estimates rather than relying on adaptive scaling alone,
both versions address noise-limited late-stage accuracy while
preserving the accelerated adaptive structure. 
A unified Lyapunov analysis handles the dependence between the
recursive estimate and its synchronous preconditioner. With the
prescribed polynomial steps, gradient-estimation errors diminish,
and the objective bounds retain accelerated initial-error decay
alongside a vanishing stochastic contribution. The same analysis
also establishes convergence of projected nonrecursive
{\it Adam-SHANG-s} under polynomial steps, filling a gap in its
original theory.

Future directions include extending the framework to nonconvex
objectives, exploring projection-free variants for greater flexibility,
and further studying the stochastic behavior of preconditioner-based
spectral steps. Incorporating the directional ratio $q_k$ into step
selection may also enable less conservative updates and sharper
convergence guarantees for coordinatewise adaptation.

\bibliographystyle{elsarticle-num}
\bibliography{references}

\appendix
\section{Details of the convex experiments}\label{app:convex-experiments}
\subsection{Tuning details}\label{app:common-implementation}
The grids for Adam, STORM, and $\mu^2$-ExtraSGD are 
\begin{align*}
\ell_0&\in\{5,1,0.5,0.1,0.05,0.01,0.005,0.001,0.0005,0.0001\},\\
c/L^2&\in\{28.00001,28.0001,28.001,28.01,28.1,29,32,64,128\},\\
\theta&\in\{0.001,0.005,0.01,0.05,0.1,0.5,1\}.
\end{align*}
For STORM, the problem-specific minibatch-variance bound $\sigma^2$ sets
$$
\kappa=\left(\frac{\sigma^2}{7L(c-28L^2)}\right)^{1/3},\qquad
w=\max\left\{(4L\kappa)^3,\,2\sigma^2,\,\left(\frac{c\kappa}{4L}\right)^3\right\}.
$$
Here $w$ is an offset, not weight decay.

All experiments use pilot seeds $100,101,102$ for tuning and
evaluation seeds $5000,\ldots,5009$ for the reported curves.
Candidates are ranked by their mean terminal objective over the
three pilot runs, evaluated exactly for logistic regression and
using an independent observation batch for denoising.
Logistic minibatches are sampled uniformly with replacement;
recursive gradient differences reuse the same samples at both queries.
Adaptive SHANG implementations share an independent initialization
minibatch and set $\eta_0=0$; its $B$ sample-gradient evaluations
are included in each method's budget. 

\subsection{Graph-signal denoising}\label{app:graph-details}
In \eqref{eq:experiment-objective}, $D$ weights data fit and
unit-weight four-neighbor edges penalize neighboring differences.
For $r,j=0,\ldots,11$, set
$$
D_{ii}=12\cdot1000^{j/11}/ \textstyle \sum_{\ell=0}^{11}1000^{\ell/11},
\qquad i=12r+j+1.
$$
For $(u,v)=(j/11,r/11)$, set
$$
b(u,v)=0.45\cos(\pi u)\cos(\pi v)+0.25\sin(2\pi u)\sin(\pi v)
+0.2\exp\left(-\frac{(u-0.7)^2+(v-0.3)^2}{0.03}\right).
$$
The graph, $D$, and $b$ are fixed across noise levels.
Larger weights mean greater observation precision:
$\operatorname{Var}(z_i-b_i)=\nu^2/D_{ii}$.
We sample Gaussian minibatch means directly as
$\bar z=b+\nu B^{-1/2}D^{-1/2}\zeta$, $\zeta\sim N(0,I)$,
and use
\begin{equation}\label{eq:experiment-sample}
\nabla f(x;\bar z)=D(x-\bar z)+L_Gx.
\end{equation}
Each evaluation counts as $B$ sample gradients.
The noise $D(b-\bar z)$ has variances
\begin{equation}\label{eq:experiment-variance}
\sigma_i^2=\nu^2D_{ii}/B,\qquad
\sigma^2=\E\|\nabla f(x;\bar z)-\nabla f(x)\|^2=\nu^2d/B.
\end{equation}
$H=D+L_G$ is the common Hessian, with
$L=\lambda_{\max}(H)\approx11.072$ and condition number $171.023$.
Exact gaps are $\frac12(x-x^\star)^\top H(x-x^\star)$, where
$Hx^\star=Db$.

\subsection{Logistic regression on synthetic data}\label{app:logistic-details}
The feature-frequency profile is
$$
\pi_j=6\cdot100^{(j-1)/59}/\textstyle \sum_{\ell=0}^{59}100^{\ell/59},\qquad j=1,\ldots,60.
$$
The imbalance concerns feature frequencies, not class proportions.
Column supports are sampled uniformly without replacement, with
counts preserving the prescribed total sparsity.
Labels follow
$$
\Pr(b_i=1\mid a_i)=\frac1{1+\exp(-a_i^\top\theta)},
$$
where $\theta$ is drawn once with independent entries in
$\{-0.25,0.25\}$ and then fixed.
The stated $L$ bounds every sample Hessian. STORM uses
$
\sigma^2=\frac1{BN}\textstyle\sum_{i=1}^N\|a_i\|^2=6/B,
$
since each sample logistic gradient has norm at most $\|a_i\|$
and the regularizer is deterministic.

\subsection{Posterior radius verification}\label{app:radius-verification}
References for evaluation and verification are computed after
optimization by a linear solve (denoising), Newton's method
(synthetic logistic regression), or L-BFGS (CIFAR-10 logistic regression);
they do not determine parameters.
Strong convexity gives
$\|x_{\mathrm{ref}}-x^\star\|\le\varepsilon_{\mathrm{ref}}$, with
$$
\varepsilon_{\mathrm{ref}}=
\begin{cases}
	\|Hx_{\mathrm{ref}}-Db\|/\lambda_{\min}(H),&\text{denoising},\\
	\|\nabla f(x_{\mathrm{ref}})\|/\lambda,&\text{logistic regression}.
\end{cases}
$$

\begin{figure}[!htbp]
	\centering
	\includegraphics[width=0.45\linewidth]{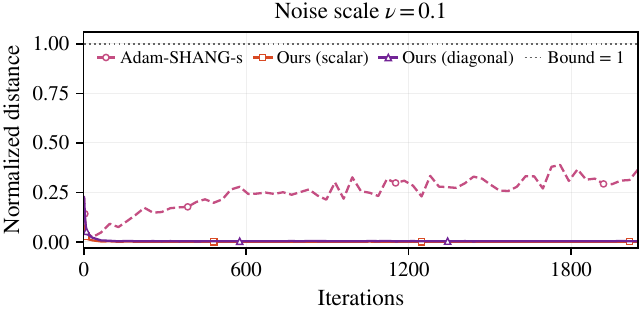}\hfill%
	\includegraphics[width=0.45\linewidth]{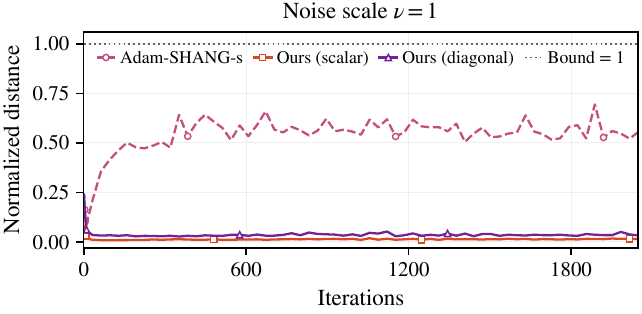}
	\par\vspace{0.4em}
	\includegraphics[width=0.45\linewidth]{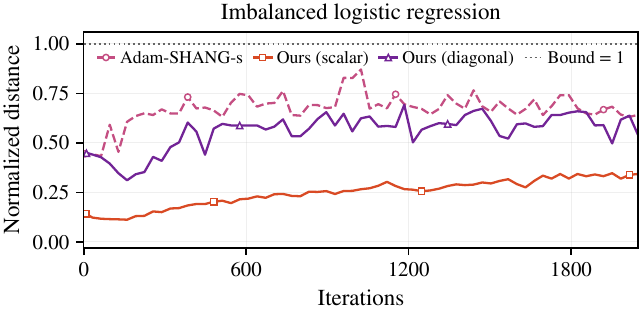}\hfill%
	\includegraphics[width=0.45\linewidth]{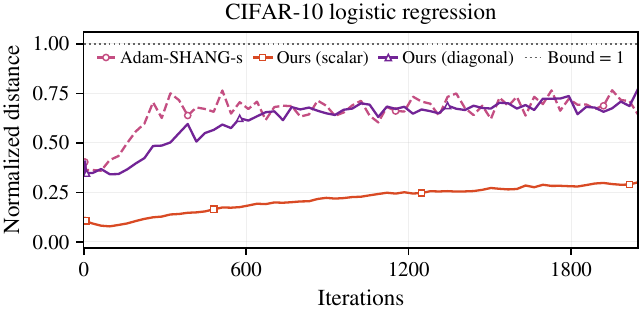}
	\caption{Posterior boundedness for graph-signal denoising (top),
		synthetic logistic regression (bottom left), and CIFAR-10 logistic regression (bottom right).
		All reported ratios remain below $1$.}
	\label{fig:boundedness}
	\vspace{-0.8em}
\end{figure}

For the prescribed radii, the reference norms (rounded below)
and residual bounds certify
$\|x_{\mathrm{ref}}\|+\varepsilon_{\mathrm{ref}}\le R/2$ (scalar) and
$\|x_{\mathrm{ref}}\|_\infty+\varepsilon_{\mathrm{ref}}\le R/2$ (diagonal).
\begin{center}
\begin{tabular}{lccc}
\hline
Problem & $\|x_{\mathrm{ref}}\|$ & $\|x_{\mathrm{ref}}\|_\infty$ & Bound on $\varepsilon_{\mathrm{ref}}$\\
\hline
Denoising & $1.950$ & $0.442$ & $5\times10^{-14}$\\
Synthetic logistic & $2.266$ & $0.888$ & $2\times10^{-13}$\\
CIFAR-10 logistic & $11.517$ & $0.674$ & $2\times10^{-5}$\\
\hline
\end{tabular}
\end{center}
Projection then gives $\|y_k-x^\star\|\le R$ (scalar) and
$\|y_k-x^\star\|_\infty\le R$ (diagonal).
Figure~\ref{fig:boundedness} plots the maximum over runs during
the first $K=2048$ iterations of
\begin{equation}\label{eq:posterior-radius-ratios}
\frac{\|y_k-x_{\mathrm{ref}}\|+\varepsilon_{\mathrm{ref}}}{R}
\quad\text{(scalar)},\qquad
\frac{\|y_k-x_{\mathrm{ref}}\|_\infty+\varepsilon_{\mathrm{ref}}}{R}
\quad\text{(diagonal)}.
\end{equation}


\end{document}